\documentclass[a4paper,11pt,reqno]{amsart}
\usepackage{amsmath}
\usepackage{amsthm}

\usepackage{graphicx}
\usepackage{amssymb}

\usepackage{color}

\usepackage[italian,english]{babel}

\usepackage[utf8x]{inputenc}
\usepackage{fancyhdr}

\usepackage{calc}

\usepackage[margin=2.6cm]{geometry}

\usepackage{srcltx}

\usepackage[T1]{fontenc}

\makeatletter
\def\cleardoublepage{\clearpage\if@twoside \ifodd\c@page\else%
		 \hbox{}%
	 \thispagestyle{empty}%              % Empty header styles
	 \newpage%
	 \if@twocolumn\hbox{}\newpage\fi\fi\fi}
\makeatother

\DeclareMathOperator{\sign}{sign}

\let\cleardoublepage\clearpage

\addto\captionsitalian{}

\newtheorem{thm}{Theorem}[section]
\newtheorem{cor}[thm]{Corollary}
\newtheorem{lem}[thm]{Lemma}
\newtheorem{pro}[thm]{Proposition}

\newtheorem{oss}[thm]{Remark}
\numberwithin{equation}{section}

\newcommand{\Hr}{H^1_{\textnormal{rad}}}

\begin{document}

\title[Radial fast diffusion on the hyperbolic space]{Radial fast diffusion on the hyperbolic space}
\author {Gabriele Grillo, Matteo Muratori}

\address {Gabriele Grillo, Matteo Muratori: Dipartimento di Matematica, Politecnico di Milano, Piaz\-za Leonardo da Vinci 32, 20133 Milano, Italy}

\email {gabriele.grillo@polimi.it}
\email {matteo.muratori@polimi.it}

%\author {Gabriele Grillo, Matteo Muratori}
%\keywords{Weighted porous media equation; weighted Sobolev inequalities; nonlinear diffusion equations; smoothing effect; asymptotics of solutions.}
%
%\address {Gabriele Grillo, Matteo Muratori: Dipartimento di Matematica, Politecnico di Milano, Piaz\-za Leonardo da Vinci 32, 20133 Milano, Italy}
%
%\email {gabriele.grillo@polimi.it}
%
%
%\email {matteo1.muratori@mail.polimi.it}
%

\maketitle

\begin{abstract}
We consider positive radial solutions to the fast diffusion equation $u_t=\Delta (u^m)$ on the hyperbolic space $\mathbb{H}^{N}$ for $N \ge 2$, $m\in\left(m_s,1\right)$, $m_s=\frac{N-2}{N+2}$. By radial we mean solutions depending only on the geodesic distance $r$ from a given point $o \in \mathbb{H}^N$. We investigate their fine asymptotics near the extinction time $T$ in terms of a separable solution of the form ${\mathcal V}(r,t)=(1-t/T)^{1/(1-m)}V^{1/m}(r)$, where $V$ is the unique positive energy solution, radial w.r.t.\ $o$, to $-\Delta V=c\,V^{1/m}$ for a suitable $c>0$, a semilinear elliptic problem thoroughly studied in \cite{MS08}, \cite{BGGV}. We show that $u$ converges to ${\mathcal V}$ in relative error, in the sense that $\left\|{u^m(\cdot,t)}/{{\mathcal V}^m(\cdot,t)}-1\right\|_\infty\to0$ as $t\to T^-$. In particular the solution is bounded above and below, near the extinction time $T$, by multiples of $(1-t/T)^{1/(1-m)}e^{-(N-1)r/m}$. Solutions are smooth, and bounds on derivatives are given as well. In particular, sharp convergence results as $t\to T^-$ are shown for spatial derivatives, again in the form of convergence in relative error.
\end{abstract}

\section{Introduction and preliminaries}\label{sec: Basics}
We analyze the asymptotic behaviour of solutions to the following \emph{Fast Diffusion Equation} (FDE) on the $N$-dimensional hyperbolic space ${\mathbb H}^N$ (throughout the whole of this work we shall always assume $N\ge2$):
\begin{equation} \label{eq: Fast}
\begin{cases}
u_t = \Delta \left( {u^m} \right)   &  \textnormal{on }   {\mathbb H}^N \times (0,T)   \\
%(u^m)_r=0   &  \textnormal{for } (r,t) \in \left\{ 0 \right\} \times (0,T)   \\
u =u_0 \geq 0   &  \textnormal{on }   {\mathbb H}^N \times \{ 0 \}
 \end{cases} \, ,
\end{equation}
where $\Delta$ is the Riemannian Laplacian, $m \in (m_s,1)$ and $m_s$ is the critical exponent
$$ m_s= \frac{N-2}{N+2} \, .  $$
The initial datum $u_0$ is assumed to be \it radial \rm in the sense that it depends only on the geodesic distance $r$ from a given point $o$ (which we shall indicate as $d(x,o)$), such a point being considered as fixed. Solutions to the FDE corresponding to radial data are of course radial as well for any fixed time. In \eqref{eq: Fast} the parameter $T=T(u_0)$ denotes, for an appropriate class of data, the \it extinction time \rm of the solution $u$, namely the smallest positive time $t$ at which $u(x,t)= 0$ identically. In fact, the results of \cite{BGV} show that, in a class of Cartan-Hadamard manifolds which includes ${\mathbb H}^N$, such a time does exist finite for initial data which belong to $L^q({\mathbb H}^N)$, where $q>\max\left(1,{N(1-m)}/2\right)$. We refer to such paper also for the relevant existence and uniqueness results provided there for solutions to \eqref{eq: Fast} (see however an alternative approximation procedure sketched in Section \ref{Sec: asy}).

It is well-known that the fine asymptotics of solutions to the FDE posed in the whole Euclidean space ${\mathbb R}^N$ is governed by (suitable rescalings of) Barenblatt, or pseudo-Barenblatt, solutions. A huge literature on the topic has been produced in the last decade, and we limit ourself to quote the recent book \cite{Vaz07} and, without any claim of completeness, the papers \cite{CT, C, DD, CL, LM, DM, CDT, MS, DS, BBDGV, BGV2, BDGV, DT, DKM}  and references quoted therein. Notice that extinction in finite time holds, in this context and for appropriate class of data, only for $m<(N-2)/N$ (provided one takes, for simplicity, $N\ge3$).

The situation on negatively curved manifolds is different from the Euclidean one since vanishing of solutions in finite time often occurs not only for $m$ close to 0, but for all $m<1$. In fact, this is instead somewhat similar to what happens in the case of the homogeneous Dirichlet problem on \it bounded \rm Euclidean domains $\Omega$, which was deeply investigated in \cite{DKV91, FS, BH80, BGVFast, SV}, where at various levels of detail it is shown that the asymptotics of suitable classes of solutions can be discussed in terms of \it separable solutions \rm of the form $(1-t/T)^{1/(1-m)}S^{1/m}(x)$, $S$ being a positive solution to the elliptic problem $-\Delta S=c\,S^{1/m}$ in $\Omega$, $S=0$ on $\partial\Omega$, for an appropriate value of the constant $c>0$ (in principle depending on the initial datum).

No positive solution to the above elliptic problem exists in the whole Euclidean space $\mathbb{R}^N$, but the fact the the bottom of the $L^2$ spectrum of $-\Delta$ on ${\mathbb H}^N$ is strictly positive points towards existence of such a solution in ${\mathbb H}^N$. In fact, this result has been proved in \cite{MS08}. More precisely, it is shown there that, given any $c>0$ and $m\in\left(m_s,1\right)$, the equation
\begin{equation}\label{emden}
-\Delta V = c \, V^{\frac{1}{m}} \ \ \ \textnormal{on $\mathbb{H}^N$}
\end{equation}
admits strictly positive solutions $V$ belonging to the energy space $H^1({\mathbb H}^N)$ (what we call \emph{energy solutions}), which are necessarily \it radial \rm with respect to some point $o$, the latter being therefore the only free parameter characterizing such solutions. Notice that solutions to \eqref{emden} associated to different values of $c$ are related by scaling, namely they are all multiples of the solution corresponding to $c=1$. The asymptotics of $V$ as $r\to \infty$ has been studied as well in \cite{MS08} and slightly improved in \cite{BGGV}: the main result is the existence of constants $l=l(c,m,N) >0$ such that, for all $k\in{\mathbb N}\cup \{0\}$,
\begin{equation}\label{eq:asympt}
\lim_{r\to \infty}e^{(N-1)r}\frac{{\rm d}^{k}V}{{\rm d}r^k}(r)=(-1)^k(N-1)^kl \, .
\end{equation}
Actually, in \cite{MS08} and \cite{BGGV}, \eqref{eq:asympt} is shown in to be valid for $k=0,1$ only, but the equation satisfied by $V$ allows to prove it for any $k \in \mathbb{N}$.

Infinitely many other positive solutions to \eqref{emden} exist, but none of them belongs to $H^1({\mathbb H}^N)$ and their behaviour as $r \to \infty$ is polynomial (see again \cite{BGGV}).

The asymptotics of a given solution $u$ to \eqref{eq: Fast} starting from a radial initial datum $u_0$ is related to the energy solution $V$ of \eqref{emden} having the same pole $o$ as $u_0$ and corresponding to a choice of $c$ that depends on $u_0$ itself via the extinction time of $u$, namely the one that satisfies
\begin{equation}\label{eq: statProb}
-\Delta V = \frac{1}{(1-m)T} \, V^{\frac{1}{m}} \ \ \ \textnormal{on $\mathbb{H}^N$}.
\end{equation}

In fact, our main result is the following Theorem \ref{ThConvRel}. Formula \eqref{eq: ConvRelFin} in it will be proved in Section \ref{Sec: rel}, while formulae \eqref{deriv1} and \eqref{eq: convder} will be proved in Section \ref{Sec: deriv}.
\begin{thm}[convergence in relative error and convergence of derivatives]\label{ThConvRel}
Let $u$ be the solution to the fast diffusion equation \eqref{eq: Fast} corresponding to a non-identically zero initial datum $u_0\ge0$, which is radial w.r.t.\ $o\in{\mathbb H}^N$ and belongs to $L^q({\mathbb H}^N)$ for some $q>{N(1-m)}/2$ with $q\ge1$. If $T>0$ is the extinction time of $u$ and $V$ is the unique positive energy solution, with pole $o$, to the stationary elliptic problem \eqref{eq: statProb}, then
\begin{equation}\label{eq: ConvRelFin}
\lim_{t \rightarrow T^{-}} \left\| \frac{u(t)}{\left( 1-\frac{t}{T} \right)^{\frac{1}{1-m}} V^{\frac1m}} - 1 \right\|_\infty = 0 \, .
\end{equation}
Moreover, for all $k\in{\mathbb N}$ there holds
\begin{equation}\label{deriv1}
\lim_{t \rightarrow T^{-}}\left|\frac{u(t)}{\left( 1-\frac{t}{T} \right)^{\frac{1}{1-m}} V^{\frac1m}}\right|_{C^k(\mathbb{R})}=0 \, ,
\end{equation}
where it is understood that
$$\left| \varphi \right|_{C^k(\mathbb{R})}:=\left\| \frac{\partial^k \varphi}{\partial r^k} \right\|_\infty$$
for all regular functions $\varphi$. As a consequence, for all $k\in{\mathbb N}$ there exists a smooth function $F_k(r)$, having the property
\[\lim_{r \to \infty} F_k(r) = 1 \, , \]
such that
\begin{equation}\label{eq: convder}
\lim_{t \rightarrow T^{-}}\left\|\frac{\frac{\partial^k u(t)}{\partial r^k}}{(-1)^k \left(\frac{N-1}{m} \right)^k\left( 1-\frac{t}{T} \right)^{\frac{1}{1-m}} V^{\frac1m}}-F_k\right\|_\infty=0 \, .
\end{equation}
\end{thm}
\begin{oss}\rm
Formulas \eqref{eq: convder}, \eqref{eq:asympt} (and identity \eqref{eq: potenzeV} below) imply that, given $k \in \mathbb{N} $, for all $\epsilon>0$ there exist $t_\epsilon\in(0,T)$ and $r_\epsilon>0$ such that
\begin{equation}\label{eq: largert}
1-\epsilon\le\frac{\frac{\partial^k u}{\partial r^k}(r,t)}{\left( 1-\frac{t}{T} \right)^{\frac{1}{1-m}}  \frac{\mathrm{d}^k V^{\frac1m}}{\mathrm{d}r^k}(r)}\le1+\epsilon\ \ \ \ \forall r\ge r_\epsilon \, , \ \forall t\in[t_\epsilon, T) \, .
\end{equation}
Notice that \eqref{eq: convder} bears some similarity with some of the results given, for the Euclidean case and in the range of $m$ for which there is no extinction, in \cite{LV}. See also \cite{LPV} for similar results for the Euclidean $p$-Laplacian driven evolution.
\end{oss}
%\begin{equation}
%\lim_{t \rightarrow T^{-}}\left|
%
%\frac{\frac{\partial {u(t)}}{\left( 1-\frac{t}{T} \right)^{\frac{m}{1-m}} V}
%
%\right|
%\end{equation}

The method proof of Theorem \ref{ThConvRel} and the known behaviour at infinity of $V$ and its derivatives allow to show that the next Theorem \ref{CorHarnack} holds. In it, we shall first state a \it global Harnack principle\rm, in the spirit of \cite{DKV91} and \cite{BV}. Secondly, we shall give upper and lower bounds on derivatives of the solution. In fact, \eqref{eq: ghp} will follow from Propositions \ref{thm: stimaSopraFond} and \ref{thm: below}, \eqref{eq: derivSP}-\eqref{eq: derivTE} follow from the results given in Section \ref{Sec: deriv} whereas \eqref{eq: lower} is an immediate consequence of \eqref{eq: largert}, \eqref{eq:asympt} and \eqref{eq: potenzeV}.

\begin{thm}[global Harnack principle and bounds for derivatives]\label{CorHarnack}
Let the assumptions of Theorem \ref{ThConvRel} be valid. Then for all $\varepsilon>0$ there exist positive constants $c_1=c_1(u_0,m,N,\varepsilon)$, $c_2=c_2(u_0,m,N,\varepsilon)$ such that the bound
\begin{equation}\label{eq: ghp}
c_1\left( 1-\frac{t}{T} \right)^{\frac{1}{1-m}}e^{-\frac{N-1}m r}\le u(r,t)\le c_2\left( 1-\frac{t}{T} \right)^{\frac{1}{1-m}}e^{-\frac{N-1}mr}
\end{equation}
holds true for all $r\ge0$, $t\in[\varepsilon,T)$. Moreover, for all $k\in{\mathbb N}$ there hold
\begin{gather}
\left|\frac{\partial^k}{ \partial r^k}u(r,t)\right| \le C_{1,k} \left( 1-\frac{t}{T} \right)^{\frac{1}{1-m}}e^{-\frac{N-1}m r} \label{eq: derivSP} \, ,  \\
\left|\frac{\partial^k}{ \partial t^k}u(r,t) \right|   \le  C_{2,k} \, \frac{e^{k(N-1)\left(\frac1m-1 \right) r } }{(1-\frac{t}{T})^k} \, e^{-\frac{N-1}{m} r} \left( 1 - \frac{t}{T} \right)^{\frac{1}{1-m}} \label{eq: derivTE}   \\
\forall r \ge0 \, , \ \forall t\in[\varepsilon,T) \nonumber
\end{gather}
for suitable positive constants $C_{1,k}=C_{1,k}(u_0,m,N,\varepsilon)$, $C_{2,k}=C_{2,k}(u_0,m,N,\varepsilon)$. In addition, for any $k\in{\mathbb N}$ there exist $\bar{t} \in(0,T)$, $\bar{r}>0$ and a suitable positive constant $C^\prime_k=C^\prime_k(u_0,m,N)$ such that
\begin{equation}\label{eq: lower}
\left|\frac{\partial^k}{ \partial r^k}u(r,t)\right|\ge  C^\prime_k \left( 1-\frac{t}{T} \right)^{\frac{1}{1-m}}e^{-\frac{N-1}m r} \ \ \ \forall r \ge \bar{r} \, , \ \forall t \in \left[ \bar{t},T \right) .
\end{equation}
\end{thm}
Notice that our global Harnack principle (as well as the corresponding bounds for derivatives) is in the spirit of the one proved in the fundamental paper \cite{DKV91} by DiBenedetto, Kwong and Vespri for bounded domains of $\mathbb{R}^N$, and of the corresponding results proved by Bonforte and V\'azquez \cite{BV} for the fast diffusion equation in ${\mathbb R}^N$: in this latter case solutions can in fact be bounded above and below for all times by suitable Barenblatt solutions. Convergence in relative error was first discussed in \cite{Vaz03} for solutions to the fast diffusion equation in ${\mathbb R}^N$, the attractor in that case being still a Barenblatt solution. Later, Bonforte, Grillo and V\'azquez showed in \cite{BGVFast} that convergence in relative error to a separable solution occurs in the case of bounded domains, thus improving the results of \cite{DKV91}.

It is worth pointing out that the techniques of proof of Theorems \ref{ThConvRel} and \ref{CorHarnack} can be used to capture the \it spatial \rm behaviour of solutions, for any \emph{fixed} $t>0$, to the FDE also in the \emph{subcritical} range $m\in\left(0, m_s \right]$. Indeed the following result holds (see Remarks \ref{oss: AllM}, \ref{oss: AllMBelow} and \ref{oss: mSottocritico} for a sketch of proof).

\begin{thm}[spatial behaviour for subcritical $m$]\label{CorSub}
Let the assumptions of Theorem \ref{ThConvRel} be valid, and suppose that $m$ lies in the subcritical range $\left(0,m_s\right]$. Then for any fixed $t\in(0,T)$ there exist positive constants $c_1=c_1(t,u_0,m,N)$, $c_2=c_2(t,u_0,m,N)$ such that the bound
\begin{equation}\label{eq: sub}
c_1 \,e^{-\frac{N-1}m r}\le u(r,t)\le c_2 \,e^{-\frac{N-1}m r}
\end{equation}
holds true for all $r\ge0$. Moreover, for all $k\in{\mathbb N}$ the bounds
\begin{gather}
\left|\frac{\partial^k}{\partial r^k}u(r,t)\right|\le C_{1,k} \, e^{-\frac{N-1}m r} \label{eq: derivSPCri} \, , \\
\left|\frac{\partial^k}{ \partial t^k}u(r,t) \right|   \le  C_{2,k} \, e^{k(N-1)\left(\frac1m-1 \right) r }  \, e^{-\frac{N-1}{m} r}  \label{eq: derivTECri} \\
\forall r\ge0,\ \forall t \in (0,T)  \nonumber
\end{gather}
hold true for suitable positive constants $C_{1,k}=C_{1,k}(t,u_0,m,N)$, $C_{2,k}=C_{2,k}(t,u_0,m,N)$.
\end{thm}

Some words have to be said about the assumption of radiality that we require on the initial data, which is related to technical issues. Firstly, we need some a priori decay properties for the solution in order to exploit suitable \emph{barrier arguments}. Such decay properties hold automatically for radial functions in the energy space but need not be valid for general solutions. In second place, it is not obvious that the solution (suitably rescaled in time) corresponding to a nonradial datum selects a unique limiting spatial profile $V$ along subsequences (recall the degree of freedom given by the pole $o$). This was proved in \cite{FS} in the Euclidean case (on bounded domains) but it is not known in the present context. Besides, in the proof of the key  Lemma \ref{lem: convLoc}  compactness of the embedding $H^1_{{\rm rad}}({\mathbb H}^N) \hookrightarrow  L^{\frac{m+1}{m}}_{{\rm rad}}({\mathbb H}^N)$ is used in a crucial way, and such property fails in the nonradial case.

Finally, notice that it is not even clear how to consider data which are not radial but bounded above and below by suitable radial data, since the extinction times of the corresponding solutions in principle change. The existence of \it ordered \rm radial data such that the corresponding solutions have the \it same \rm extinction time $T$ is an open problem. Should such a construction be possible, the methods of the present paper would give convergence in relative error to the separable solution extinguishing at time $T$ also for nonradial data in between.

\begin{oss}\rm
Keeping the fundamental hypothesis of radiality, our results hold in somewhat more general geometric frameworks, but we preferred to state them in the case of ${\mathbb H}^N$ to avoid bothering the reader with heavier notation and technicalities. In fact one could consider \it Riemannian models \rm (see \cite{GW, Besse} as general references and \cite{BFG} for the analysis of Lame-Emden-Fowler equations in such context) whose metric is defined, in spherical coordinates about a pole $o$, by $ds^2=dr^2+\psi^2(r)d\Theta^2, \, \Theta\in {\mathbb S}^{N-1}$, where $\psi\in C^2([0,\infty))$, $\psi(0)=\psi^{\prime\prime}(0)=0$, $\psi'(0)=1$, $\psi'(r)> 0$ for every $r>0$ and $\lim_{r\rightarrow \infty}\psi'(r)/\psi(r)\in(0,\infty)$. Notice that sectional curvatures at a point $P$ tend, as the geodesic distance $d(o,P)$ tends to $\infty$, to a strictly positive constant. In such a kind of manifold a radial energy solution having the properties of the present solution $V$ has been shown to exist in \cite{BFG}.
\end{oss}

\subsection{Preliminaries} \label{subs: pre}
As for the initial datum $u_0=u_0(r)$ in \eqref{eq: Fast}, in principle besides its nonnegativity we should also assume that it is bounded and such that $u_0^m \in \Hr({\mathbb H}^N)$, where
\begin{equation}\label{eq: asmpEnergy}
\Hr({\mathbb H}^N) \! = \! \left\{  v \ \textrm{radial}\! :  \| v \|_{H^1}^2 =  \int_0^\infty { v^2(s) \, (\sinh s)^{N-1} \mathrm{d}s }  + \int_0^\infty {[v^\prime(s)]^2   \, (\sinh s)^{N-1} \mathrm{d}s}  < \infty  \right\}  .
\end{equation}
Notice that $\Hr({\mathbb H}^N)$ coincides with the space of radial functions about $o$ which belong to $H^1(\mathbb{H}^N)$. By \emph{energy solutions} to \eqref{eq: Fast} one should mean those starting from data $u_0$ as in \eqref{eq: asmpEnergy} (in the nonradial context, those starting from $u_0: u_0^m \in H^1(\mathbb{H}^N) $), but in fact the results of \cite{BGV} show that the solution $u$ corresponding to an initial datum which fulfils the integrability conditions of Theorem \ref{ThConvRel} automatically satisfies $u^m(\cdot,\varepsilon) \in  \Hr({\mathbb H}^N) \cap L^\infty({\mathbb H}^N) $ for all $\varepsilon>0$. This is stated in \cite{BGV} for $N\ge3$, but it holds true when $N=2$ as well because the methods of proof exploited in \cite{BGV} rely only on the validity of a suitable Sobolev inequality in $H^1({\mathbb H}^N)$, which is valid also when $N=2$.

Let us see now what problem \eqref{eq: Fast} reads like for radial solutions. Recall that the Riemannian Laplacian on the hyperbolic space, for radial functions $v=v(r)$, takes the form
\begin{equation} \label{eq: RadLap}
\Delta {v(r)}= \frac1{(\sinh r)^{N-1}} \,  \left[ (\sinh r)^{N-1} \, v^\prime(r)  \right]^\prime = v^{\prime\prime}(r) + (N-1)(\coth r)v^\prime(r) \, ,
\end{equation}
where the apex $^\prime$ stands for derivation w.r.t.\ $r$. From \eqref{eq: RadLap} we have that studying  energy solutions to \eqref{eq: Fast} for radial initial data is equivalent to studying energy solutions to the problem
\begin{equation} \label{eq: FastRadial}
\begin{cases}
u_t = \left(u^m\right)^{\prime\prime} + (N-1)(\coth r) \left(u^m\right)^\prime &  \textnormal{ in }  (0,\infty) \times (0,T)   \\
(u^m)^\prime =0   &  \textnormal{ on }  \{ 0 \} \times (0,T)   \\
u =u_0 \geq 0   &  \textnormal{ on }  [0,\infty)  \times \{ 0 \}
\end{cases} \, .
\end{equation}

The fact that there exists a \emph{finite} extinction time $T>0$ is a straightforward consequence of the validity in $\Hr(\mathbb{H}^N)$ (indeed also in $H^1(\mathbb{H}^N)$) of both a Poincar\'e and a Sobolev inequality (see, for instance, \cite[Sect. 5.10]{Vaz07} and \cite[Sect. 3]{BGGV}), that is
\begin{equation}\label{eq: poinSob}
\| v \|_{2} \le C_P \, \| v^\prime \|_{2}  \, , \ \ \ \| v \|_{\frac{2N}{N-2}} \le C_S \, \| v^\prime \|_{2}
\end{equation}
for all $v \in \Hr({\mathbb H}^N)$ and suitable positive constants $C_P=C_P(N)$, $C_S=C_S(N)$, where
$$  \|v\|_{p}^p=\int_0^\infty { |v|^p(s) \, (\sinh s)^{N-1} \mathrm{d}s } \, , \ \ \ L^p_{\textnormal{rad}}({\mathbb H}^N) = \left\{ v \ \textrm{radial}\! : \, \| v \|_{p} < \infty \right\} . $$
Moreover, one can prove \cite[Th. 3.1]{BS12} that the embedding of $\Hr({\mathbb H}^N)$ into $L^p_\textnormal{rad}(\mathbb{H}^N)$ is \emph{compact} for all $p \in \left(2,2N/(N-2)  \right)$. Notice that, since $m \in (m_s,1)$, this means in particular that
\begin{equation}\label{eq: compEmb}
\Hr({\mathbb H}^N) \Subset L^{\frac{m+1}{m}}_{\textnormal{rad}}({\mathbb H}^N) \, ,
\end{equation}
a crucial fact that we shall exploit in the next section. Recall however that the compact embedding \eqref{eq: compEmb} \emph{fails} in $H^1(\mathbb{H}^N)$, another nontrivial issue that points out the advantage of working in the radial framework.

In the sequel, for notational simplicity, we shall write $L^p$ instead of $L^p({\mathbb H}^N)$ and do the same for all the functional spaces involved.
\subsection{Plan of the paper} All the above results will be proved in several intermediate steps. Local uniform convergence of $u^m/(1-t/T)^{m/(1-m)}$ to the stationary solution $V$ is shown in Section \ref{Sec: asy}, along lines similar to the ones \cite{BH80}. Then, a suitable upper bound for solutions (holding as $ r \to \infty$) is shown in Section \ref{Sec: Above}, whereas a matching lower bound is proved in Section \ref{sec: below}. The more delicate issue, namely the passage to the relative error $u^m/[(1-t/T)^{m/(1-m)} V] - 1 $, is dealt with in Section \ref{Sec: rel}. Section \ref{Sec: deriv} contains the proofs of the results concerning space-time derivatives of solutions, which exploit both regularity theory and the claimed convergence in relative error.

\section{Local uniform convergence of the rescaled solution to the stationary profile}\label{Sec: asy}
As previously mentioned, each solution to \eqref{eq: Fast} extinguishes in a finite time $T>0$. Therefore the asymptotic behaviour of $u$ is, from this point of view, trivial: the solution goes to zero as $t\uparrow T$. In order to study finer properties of $u$ it is very useful to look for \emph{separable solutions} to \eqref{eq: Fast} (if any), so that their asymptotic behaviour might unveil at least the expected order of convergence to zero of a generic solution. Hence, let us set $u(x,t)=g(t) \, V^{1/m}(x)$. After some straightforward computations one gets that $u$ is a solution to \eqref{eq: Fast} for some $u_0 \ge 0$ (not identically zero) if and only if $g$ satisfies
\begin{equation}\label{eq: AsTemp}
g(t)=\left( 1 - \frac{t}{T} \right)^{\frac{1}{1-m}}  \ \ \ \forall t \in [0,T]
\end{equation}
and $V$ is a positive solution to the elliptic problem \eqref{eq: statProb} for some parameter $T>0$ (the extinction time). When $m \in (m_s,1)$ existence and uniqueness of such a $V$ and its dependence the sole radial coordinate $r$ is guaranteed by compactness and by a moving plane method (see the fundamental paper \cite{MS08}). Local regularity and strict positivity of $V$ are instead a consequence of standard elliptic arguments. So the velocity of convergence to zero as $t \uparrow T$ for separable solutions is given by \eqref{eq: AsTemp}. This suggests that, in order to analyze a nontrivial asymptotics, it is convenient to study the behaviour of the rescaled solution $u(r,t)/g(t)$. Notice that, if $u(r,t)=g(t) \, V^{1/m}(r)$, then such rescaled solution trivially coincides with $V^{1/m}$. For a generic $u$ this is of course not true: however, $V^{1/m}$ seems to naturally maintain the role of an attractor for $u/g$.

Motivated by the discussion above, given the extinction time $T$ associated to the solution $u$ of \eqref{eq: Fast}, let us consider the rescaled solution $w$ defined as
\begin{equation}\label{eq: rescSol}
w(r,\tau)=  \left( \frac{T}{T-t} \right)^{\frac{1}{1-m}} \! u(r,t) \, = \, e^{\frac{\tau}{(1-m)T}} \, u\left(r, T-Te^{-\frac{\tau}{T}} \right) , \ \ \ \tau=T \log{\left( \frac{T}{T-t} \right)}
\end{equation}
$$ \forall r \in (0,\infty)  \, , \ \forall t \in (0,T) \, , \ \forall \tau \in (0,\infty) \, .$$
Straightforward computations show that $w$ solves the following problem:
\begin{equation} \label{eq: FastRisc}
\begin{cases}
w_\tau = \Delta \left( {w^m} \right) + \frac{1}{(1-m)T} \, w   &  \textnormal{ in } (0,\infty) \times (0,\infty)   \\
\left(w^m\right)^\prime = 0   &  \textnormal{ on } \left\{ 0 \right\}  \times (0,\infty) \\
w =u_0 \geq 0   &  \textnormal{ on }  [0,\infty) \times \left\{ 0 \right\}
\end{cases} \, .
\end{equation}
The aim of this section is to prove that $w^m$ converges \emph{locally uniformly} in $\{r \in [0,\infty)\}$ to $V$ as $\tau \to \infty$ (since $V$ is positive, this is equivalent to claiming that $w$ converges locally uniformly to $V^{1/m}$). The basic estimates one needs to exploit in order to prove such result were obtained in a celebrated paper \cite{BH80} by Berryman and Holland, though for regular solutions to the FDE on regular bounded domains of $\mathbb{R}^N $. Here, first we shall only point out how their techniques, with minor modifications, can be applied to this framework too. This will ensure local uniform convergence at least away from $ \{r=0\} $, while some further work will be required to extend the result to neighbourhoods of the origin $ o $.

To this end, it is convenient to see $u$ as a monotone increasing limit of the sequence of solutions $\{u_n\}$ (with extinction times $ \{T_n\} $) to the problems
\begin{equation} \label{eq: Fastn}
\begin{cases}
{(u_n)}_t = \Delta \left( {u_n^m} \right)   &  \textnormal{ in }  (0,n) \times (0,T_n)   \\
u_n=0   &  \textnormal{ on }  \left\{ n \right\} \times (0,T_n)  \\
(u^m_n)^\prime=0   &  \textnormal{ on } \left\{ 0 \right\} \times  (0,T_n) \\
u_n =u_{0n} \geq 0   &  \textnormal{ on } [0,n] \times \left\{ 0 \right\}
 \end{cases} \, ,
\end{equation}
where $\{u_{0n} \}$ is a sequence of regular data such that $u_{0n}(n)=0$, $(u^m_{0n})^\prime(0)=0$, $u_{0n} \le u_0 $, which suitably approximates $u_0$, and $ T_n \uparrow T$. We shall identify $u_n(\cdot,t)$ as functions in the whole $[0,\infty)$ by extending them to be zero outside $[0,n]$. Notice that \eqref{eq: Fastn} corresponds to the radial FDE with homogeneous Dirichlet boundary conditions posed on the ball of radius $n$ of $\mathbb{H}^N$ centered at $x=o$.
\begin{lem}\label{lem: stimeNorma}
There exists a positive constant $C=C(m,N)$ such that
\begin{equation}\label{eq: stimeNorma0}
  C  \left( T-t \right)^{\frac{1}{1-m}}  \le  \left\| u(t)  \right\|_{m+1}  \le  \left(1-\frac{t}{T} \right)^{\frac{1}{1-m}} \, \left\| u_0  \right\|_{m+1}   \ \ \ \forall t \in (0,T)  \, .
\end{equation}
Moreover, the ratio
\begin{equation}\label{eq: stimeNorma1}
\frac{\left\| (u^m)^\prime(t) \right\|_{2}}{\left\| u(t)  \right\|^m_{m+1,}}
\end{equation}
is nonincreasing along the evolution.
\begin{proof}
The left inequality in \eqref{eq: stimeNorma0} can be proved exactly as in \cite[Lemma 1]{BH80} using the identity 
\begin{equation}\label{deriv}
\frac{{\rm d}}{{\rm d}t}\int_0^\infty u^{m+1}(s,t) \, (\sinh s)^{N-1} \mathrm{d}s=-(m+1)\int_0^\infty [(u^m)^\prime(s,t)]^2 \, (\sinh s)^{N-1} \mathrm{d}s
\end{equation}
and the Poincar\'e-Sobolev inequalities in \eqref{eq: poinSob}.
To justify the the other statements we proceed as in \cite[Lemma 2]{BH80}, outlining only the main steps. At a formal level, we have:
\begin{equation}\label{deriv-first}
\frac{\int_0^\infty [(u^m)^\prime(s,t)]^2 \, (\sinh s)^{N-1} \mathrm{d}s}{\int_0^\infty u^{m+1}(s,t) \, (\sinh s)^{N-1} \mathrm{d}s}\le\frac{\int_0^\infty u^{m-1}(s,t) \, \left[ \Delta (u^m)(s,t)\right]^2 \, (\sinh s)^{N-1} \mathrm{d}s}{\int_0^\infty  [(u^m)^\prime(s,t)]^2\, (\sinh s)^{N-1} \mathrm{d}s}\,,
\end{equation}
\begin{equation}\label{deriv-second}
\frac{{\rm d}}{{\rm d}t}\int_0^\infty [(u^m)^\prime(s,t)]^2 \, (\sinh s)^{N-1} \mathrm{d}s=-2m\,\frac{\int_0^\infty u^{m-1}(s,t) \, \left[ \Delta (u^m)(s,t)\right]^2 \, (\sinh s)^{N-1} \mathrm{d}s}{\int_0^\infty  [(u^m)^\prime(s,t)]^2\, (\sinh s)^{N-1} \mathrm{d}s}\,,
\end{equation}
where \eqref{deriv-first} follows from integration by parts and Cauchy-Schwarz. From \eqref{deriv}, \eqref{deriv-first} and \eqref{deriv-second} one shows easily, exactly as in \cite[Lemma 2]{BH80}, that the ratio in \eqref{eq: stimeNorma1} is nonincreasing. Thanks to this, the right inequality in \eqref{eq: stimeNorma0} also follows as in \cite[Lemma 2, formula (14)]{BH80}. 

To justify such steps it is convenient to pass through the approximating solutions $ \{u_n\} $. Indeed the proof of \cite[Lemma 2]{BH80} requires the finiteness of the quantity
\begin{equation*}
\int_0^\infty u^{m-1}(s,t) \, \left[ \Delta (u^m)(s,t)\right]^2 \, (\sinh s)^{N-1} \mathrm{d}s
\end{equation*}
for all $t \in (0,T) $, which a priori may not hold here. However, one obtains \eqref{eq: stimeNorma0} and \eqref{eq: stimeNorma1} for $u_n$ (we can assume that $u_n^m$ is regular enough up to the boundary) and then passes to the limit as $n \to \infty$. This is feasible since $\{ u_n^m( t) \} $ converges weakly in $\Hr$ to $u^m (t)$, and by monotonicity $\{ u_n (t)\} $ converges to $u(t)$ in $L^{m+1}_{\textnormal{rad}}$ and $ \{T_n\} $ converges to $T$.
\end{proof}
\end{lem}
The next result is a key one in order to establish the mentioned convergence of $w^m$ to the stationary profile $V$.
\begin{lem}\label{lem: keyLem}
The following inequality holds true for all $\tau \in (0,\infty)$:
\begin{equation}\label{eq: keyIne}
\begin{split}
& \int_0^\infty \left[\frac12 \left[(w^m)^\prime(s,\tau)\right]^2 - \frac{m}{(1-m^2)T} w^{m+1}(s,\tau)  \right] \, (\sinh s)^{N-1} \mathrm{d}s  \\
& + m \int_0^\tau \! \int_0^\infty w^{m-1}(s,\sigma) \, \left[w_{\tau}(s,\sigma)\right]^2 \, (\sinh s)^{N-1} \mathrm{d}s  \, \mathrm{d}\sigma \\
\le & \int_0^\infty \left[\frac12 \left[(u_0^m)^\prime (s)\right]^2 - \frac{m}{(1-m^2)T} u_0^{m+1}(s)  \right] \, (\sinh s)^{N-1} \mathrm{d}s  \, .
\end{split}
\end{equation}
\begin{proof}
For the smooth rescaled solutions $w_n$ inequality \eqref{eq: keyIne} is in fact an equality, since by straightforward computations one verifies that
\begin{equation}\label{eq: deriva0}
\begin{split}
 & \frac{\mathrm{d}}{\mathrm{d}\tau} \int_0^\infty \left[\frac12 \left[(w_n^m)^\prime(s,\tau)\right]^2 - \frac{m}{(1-m^2)T_n} w_n^{m+1}(s,\tau)  \right] \, (\sinh s)^{N-1} \mathrm{d}s \\
 = & - m  \int_0^\infty w_n^{m-1}(s,\tau) \, \left[(w_n)_{\tau}(s,\tau)\right]^2 \, \chi_{(0,n)}(s) \, (\sinh s)^{N-1} \mathrm{d}s \, .
\end{split}
\end{equation}
To get estimate \eqref{eq: keyIne} it suffices to integrate \eqref{eq: deriva0} from $0$ to $\tau$ and let $n \to \infty$: to the first integral on the l.h.s. of \eqref{eq: keyIne} we can apply the weak convergence of $ \{u_n^m(t)\} $ to $ u^m(t) $ in $\Hr$ and the strong convergence of $ \{u_n(t)\} $ to $u(t)$ in $L^{m+1}_{\textnormal{rad}}$, while the second integral is handled by means of Fatou's Lemma (thanks to local regularity we can assume that, up to subsequences, $ \{(u_n)_t\} $ converges pointwise to $u_t$) or by the fact that
$$ \left\{ u_n^{\frac{m+1}{2}} \right\} \rightharpoonup u^{\frac{m+1}{2}} \ \ \ \textnormal{in } H^1(0,T; L^2_{\textnormal{rad}}) \, . $$
\end{proof}
\end{lem}
Now we are able to prove the following important result.
\begin{lem}\label{lem: convLoc}
Let $w$ be the rescaled solution \eqref{eq: rescSol} to \eqref{eq: FastRisc} and $V$ the radial, positive energy solution to the stationary problem \eqref{eq: statProb}. Then
\begin{equation}\label{eq: convLoc}
\lim_{\tau \to \infty}  \left\|  w^m(\tau)-V \right\|_{L^\infty_{\textnormal{loc}}((0,\infty))} = 0 \, ,
\end{equation}
that is $w^m(\tau)$ converges uniformly to $V$ in any compact set $ K \Subset (0,\infty)$ as $\tau \to \infty$.
\begin{proof}
We adapt the proof of \cite[Th. 2]{BH80}. First of all notice that, from \eqref{eq: keyIne}, one deduces the existence of a sequence $	\{\tau_n \} \to \infty$ such that
\begin{equation}\label{eq: provaThFon0}
\int_0^\infty w^{m-1}(s,\tau_n) \, \left[w_{\tau}(s,\tau_n)\right]^2 \, (\sinh s)^{N-1} \mathrm{d}s \rightarrow 0 \, .
\end{equation}
To show this fact notice that, thanks to \eqref{eq: stimeNorma0} and \eqref{eq: rescSol}, $\|w(\tau)\|_{m+1}$ is bounded as a function of $\tau$, hence the first integral on the l.h.s. of \eqref{eq: keyIne} is bounded from below. Moreover, the r.h.s. does not depend on $\tau$, therefore the integral
$$  \int_0^\infty \! \int_0^\infty w^{m-1}(s,\sigma) \, \left[w_{\tau}(s,\sigma)\right]^2 \, (\sinh s)^{N-1} \mathrm{d}s  \, \mathrm{d}\sigma $$
must be finite. Still from \eqref{eq: keyIne} and \eqref{eq: stimeNorma0} one gets the boundedness of $\| w^m(\tau_n) \|_{H^1}$; hence, up to subsequences, $ \{ w^m(\tau_n) \} $ converges weakly in $\Hr$ to a certain function $R$. From the compact embedding \eqref{eq: compEmb}, such convergence is in fact \emph{strong} in $L^{(m+1)/m}_{\textnormal{rad}}$. In particular, $R$ is a nonnegative non-identically zero function (indeed \eqref{eq: stimeNorma0} prevents $\| w^m(\tau_n) \|_{(m+1)/m}$ from going to zero) belonging to $\Hr$.

The next step is to show that $R$ solves \eqref{eq: statProb}. To this end, take any test function $\phi:[0,\infty) \rightarrow \mathbb{R}$ with compact support in $[0,\infty)$, multiply by it the first equation in \eqref{eq: FastRisc} (evaluated at $\tau=\tau_n$) and integrate by parts in $[0,\infty)$. This leads to the identity
\begin{equation}\label{eq: provaThFon1}
\begin{split}
\int_0^\infty w_\tau(s,\tau_n) \, \phi(s) \, (\sinh s)^{N-1} \mathrm{d}s  = &- \int_0^\infty (w^m)^\prime(s,\tau_n)  \, \phi^\prime(s) \, (\sinh s)^{N-1} \mathrm{d}s + \\
&+ \int_0^\infty  \frac{1}{(1-m)T} \, w(s,\tau_n) \, \phi(s) \, (\sinh s)^{N-1} \mathrm{d}s \, .
\end{split}
\end{equation}
The two integrals on the r.h.s. of \eqref{eq: provaThFon1} are stable under passage to the limit as $n\to \infty$: indeed $ \{w^m(\tau_n)\} $ converges weakly in $\Hr$ to $R$ and $ \{w(\tau_n)\} $ converges strongly in $L^{m+1}_{\textnormal{rad}}$ to $R^{1/m}$ (and so also in $L^1_{\textnormal{rad}}$ locally). Finally, the left hand side goes to zero since its modulus is bounded by
$$ \left( \int_0^\infty  w^{1-m}(s,\tau_n) \, \phi^2(s) \, (\sinh s)^{N-1} \mathrm{d}s \right)^{\frac12}  \left( \int_0^\infty  w^{m-1}(s,\tau_n) \, \left[w_\tau(s,\tau_n)\right]^2 \, (\sinh s)^{N-1} \mathrm{d}s \right)^{ \frac12 } \, , $$
which goes to zero thanks to \eqref{eq: provaThFon0} and to the boundedness of $\| w(\tau_n) \|_{m+1}$. From the arbitrariness of $\phi$ and the uniqueness of energy solutions to \eqref{eq: statProb} we infer that $R$ must coincide with $V$. Moreover, convergence of $ \{w^m(\tau_n)\} $ to $R=V$ in $\Hr$ is also strong. To prove that, just replace $\phi(\cdot)$ by $w^m(\cdot,\tau_n)$ in the computations above and get
\begin{equation*}\label{eq: provaThFon2}
\begin{split}
\lim_{n \to \infty}  \int_0^\infty \left[ (w^m)^\prime(s,\tau_n) \right]^2  \, (\sinh s)^{N-1} \mathrm{d}s &= \int_0^\infty  \frac{1}{(1-m)T} \, V^{\frac{m+1}{m}}(s) \, (\sinh s)^{N-1} \mathrm{d}s \\
& = \int_0^\infty [V^\prime(s)]^2  \, (\sinh s)^{N-1} \mathrm{d}s \, .
\end{split}
\end{equation*}
Hence, weak convergence plus convergence of the norms in $\Hr$ gives the claimed strong convergence. Since $\Hr$ is continuously embedded in $L^\infty_{\textnormal{loc}}$ (see e.g.\ Lemma \ref{lem: EmbLinf} below), we have proved \eqref{eq: convLoc} along the special sequence $\{\tau_n\}$. To prove that such convergence takes place along any other subsequence, one can argue by contradiction. That is, suppose there exists a sequence $\{\tau_k\}$ such that $\{w^m(\tau_k)\} $ \emph{does not} converge strongly in $\Hr$ to $V$. By \eqref{eq: keyIne} we can assume that, up to subsequences, $ \{w^m(\tau_k)\} $ converges weakly in $\Hr$ to a certain function $Q$. Now notice that, again from \eqref{eq: stimeNorma0} and \eqref{eq: keyIne} (up to a time origin shift), both
$$\left\| w(\tau) \right\|_{m+1}^{m+1}  $$
and
$$ \frac12 \left\| \left( w^m \right)^\prime \! (\tau) \right\|_{2}^2 - \frac{m}{(1-m^2)T} \left\| w(\tau) \right\|_{m+1}^{m+1}  $$
are nonincreasing functions of $\tau$. In particular,
\begin{equation}\label{eq: provaThFon3}
\left\| Q  \right\|_{\frac{m+1}{m}} = \lim_{k\to \infty} \left\| w(\tau_k) \right\|_{m+1}^m = \lim_{ n \to \infty} \left\| w(\tau_n) \right\|_{m+1}^m   = \left\| V  \right\|_{\frac{m+1}{m}}
\end{equation}
and
\begin{equation}\label{eq: provaThFon4}
\begin{split}
\left\| Q^\prime \right\|_{2}^2  &\le \liminf_{k \to \infty} \left(  \left\| \left( w^m \right)^\prime \! (\tau_k) \right\|_{2}^2 \pm \frac{m}{(1-m^2)T} \left\| w(\tau_k) \right\|_{m+1}^{m+1}  \right) \\
& = \lim_{n\to \infty} \left\| \left( w^m \right)^\prime \!  (\tau_n)  \right\|_{2}^2 = \left\| V^\prime \right\|_{2}^2 \, .
\end{split}
\end{equation}
Since $V$ is the \emph{unique} minimizer of $\| v^\prime \|_{2}^2$ among all functions $v$ with prescribed $L^{(m+1)/m}_{\textnormal{rad}}$ norm (see \cite{MS08}), \eqref{eq: provaThFon3} and \eqref{eq: provaThFon4} necessarily imply that $Q=V$. \emph{Strong} convergence of $\{w^m(\tau_k)\} $ to $Q=V$ is then a consequence of \eqref{eq: provaThFon4}, which leads to a contradiction.
\end{proof}
\end{lem}

We are left with proving that the uniform convergence \eqref{eq: convLoc} takes place also down to $r=0 $. In order to do it, we shall use Lemma \ref{lem: convLoc} and the following two lemmas, which show how positivity and boundedness of $w$ can be extended to a neighbourhood of $o$.

\begin{lem}\label{lem: posLocZero}
%Let $w$ be the rescaled solution \eqref{eq: rescSol} to \eqref{eq: FastRisc} and $V$ the radial, positive energy solution to the stationary problem \eqref{eq: statProb}.
For any $\epsilon>0$ there exist $r_\epsilon > 0$ sufficiently small and $\tau_\epsilon > 0$ sufficiently large such that
\begin{equation} \label{eq: stimaBassoIntornoZero}
w^m(r,\tau) \ge V(0)-\epsilon \ \ \ \forall(r,\tau) \in \left[0,r_\epsilon \right] \times [\tau_\epsilon,\infty) \, .
\end{equation}
\begin{proof}
We can adapt the techniques of \cite[Lem. 6.2]{DKV91}. First of all recall that, thanks to the local uniform convergence to the stationary profile \eqref{eq: convLoc}, we have uniform boundedness away from zero in any compact set which does not contain $o$. In particular, consider a point $x_0 \in \mathbb{H}^N $ such that $ r_0=d(x_0,o) \in (0,1/2) $. For a given $ \tau_0>0$, set
\begin{equation}\label{eq: defKsub}
k=\inf_{(x,\tau): \ \tau \ge \tau_0 \, , \  x \in B_{r_0/2}(x_0)  } w(d(x,o),\tau) \, ,
\end{equation}
where $ B_{r_0/2}(x_0) $ is the hyperbolic ball of radius $r_0/2$ centered at $ x_0 $. Thanks to the observations above, $k>0$ provided $ \tau_0 $ is sufficiently large. Let us consider equation \eqref{eq: FastRisc} (more precisely, its interpretation as a differential equation on $\mathbb{H}^N$) centered at $x_0$ in place of $o$. To avoid confusion, we shall call $\rho$ the radial coordinate about such $x_0$. Upon defining
$$ \widetilde{N} = 1 + (N-1) \sup_{\rho \in (0,1)} \rho \coth(\rho) \, , $$
for any function $f=f(\rho)$ such that $f^\prime(\rho) \le 0$ we have
\begin{equation}\label{eq: hypLapEucLap}
\Delta f(\rho) =  f^{\prime\prime}(\rho) + (N-1) \coth(\rho) f^\prime(\rho) \ge f^{\prime\prime}(\rho) + \frac{\widetilde{N}-1}{\rho} f^\prime(\rho)  \ \ \ \forall \rho \in (0,1) \, ,
\end{equation}
where the term on the r.h.s. of \eqref{eq: hypLapEucLap} is the Euclidean Laplacian of $f$ associated to the ``artificial dimension'' $ \widetilde{N} $. Therefore in order to seek for a \emph{subsolution} $\psi(\rho,\tau)$ to \eqref{eq: FastRisc} centered at $x_0$ it is enough to ask (we keep denoting as $^\prime$ derivative w.r.t.\ $\rho$)
\begin{equation}\label{eq: psiSubSol}
\psi_\tau \le \left( \psi^m \right)^{\prime\prime} + \frac{\widetilde{N}-1}{\rho} \left(\psi^m\right)^\prime \, , \ \  \psi \ge 0  \, ,  \ \ \left(\psi^m\right)^\prime \le 0 \, ,
\end{equation}
as long as $\rho$ varies in $(0,1)$. The proof of Lemma 6.2 of \cite{DKV91} ensures that the function
$$ \psi(\rho,\tau) = k \frac{(1-\rho^\beta)^{\frac{2}{m}}}{\left(1 + k^{1-m} \frac{b \rho^2}{\tau-\tau_0}  \right)^{\frac{\theta}{1-m}}}  $$
satisfies \eqref{eq: psiSubSol} in the region
$$ \left\{ (\rho,\tau) \in \left(\frac{r_0}{2},1\right) \times \left(\tau_0,\tau_0+\frac{r_0^2}{4} \right) \right\}  $$
upon choosing appropriately the positive parameters $\beta=\beta(m,\widetilde{N})$, $\theta=\theta(m,\widetilde{N})$ and $b=b(m,\widetilde{N},k)$, $k$ being as in \eqref{eq: defKsub}. Let us check conditions on the parabolic boundary. For $\tau=\tau_0$ and for $\rho=1$ we have, by construction, $\psi=0$ and so trivially $\psi \le w$. For $\rho=r_0/2$ (actually for any $\rho \in  (0,1)$) there holds $\psi \le k$, from which $\psi \le w$ by definition of $k$. Hence, by comparison,
\begin{equation*}
\begin{split}
w\left(d(x,o),\tau_0+{r_0^2}/{4}\right) \! & \ge \! \psi(d(x,x_0),\tau_0+{r_0^2}/{4}) \!  \\
 & \ge \! \psi(3r_0/2,\tau_0+{r_0^2}/{4}) \! = \! k \frac{\left( 1-\left(\frac{3r_0}{2}\right)^\beta \right)^{\frac{2}{m}}}{\left(1 + 9 b k^{1-m} \right)^{\frac{\theta}{1-m}}} \! = \! C_0 \! > \! 0
\end{split}
\end{equation*}
$$ \forall x \in \mathbb{H}^N: \ \ \frac{r_0}{2} \le d(x,x_0) \le \frac{3r_0}{2} \, . $$
In particular we obtain the existence of a radius $r_1>0$ and a time $\tau_1>0$ such that
\begin{equation}\label{eq: posLocBase}
w(r,\tau) \ge C_0 > 0 \ \ \ \forall (r,\tau) \in \left[0 , r_1\right] \times \{ \tau_ 1\}  .
\end{equation}
Indeed \eqref{eq: posLocBase} holds for \emph{all} $\tau$ greater than $\tau_1$, rather than only for $\tau=\tau_1$. This is a trivial consequence of the fact that $C_0$ is a subsolution to \eqref{eq: FastRisc} (the comparison condition on the lateral boundary $ \{ r_1 \} \times (\tau_1,\infty)$ is satisfied provided $C_0$ is small enough, again as a consequence of the local uniform convergence \eqref{eq: convLoc}).

Finally, we need to refine estimate \eqref{eq: posLocBase}. To this end, just observe that the function
$$ g(\tau)=C_0 \,  e^{\frac{\tau-\tau_\ast}{(1-m)T}} $$
is a solution to the differential equation in \eqref{eq: FastRisc} for any $ \tau_\ast>0 $. Still from the local uniform convergence \eqref{eq: convLoc} (and from the fact that $V(\cdot)$ is decreasing) we have that, for any $\epsilon>0$, we can choose $r_2=r_2(\epsilon)<r_1$ and $\tau_2=\tau_2(\epsilon)>\tau_1$ such that $w^m(r_2,\tau) \ge V(0)-\epsilon $ for all $\tau \ge \tau_2$. Therefore $g(\tau)$, with the choice $ \tau_\ast=\tau_2 $, is a subsolution to \eqref{eq: FastRisc} in the region
$$ \{ (r,\tau) \in (0,r_2) \times (\tau_2,\tau_3) \} \, , $$
where $\tau_3$ is the time at which $g^m(\tau_3)=V(0)-\epsilon$. Since the constant $(V(0)-\epsilon)^{1/m}$ is then a subsolution in $\{ (r,\tau) \in (0,r_2) \times (\tau_3,\infty) \}$, estimate \eqref{eq: stimaBassoIntornoZero} follows.
\end{proof}
\end{lem}

Now we prove the analogue of estimate \eqref{eq: stimaBassoIntornoZero} from above.
\begin{lem}\label{lem: boundLocZero}
For any $\epsilon>0$ there exist $r_\epsilon > 0$ sufficiently small and $\tau_\epsilon > 0$ sufficiently large such that
\begin{equation} \label{eq: boundIntornoZero}
w^m(r,\tau) \le V(0) + \epsilon \ \ \ \forall (r,\tau) \in [0,r_\epsilon] \times [\tau_\epsilon,\infty) \, .
\end{equation}
\begin{proof}
Again, we shall proceed by constructing a proper \emph{supersolution} to \eqref{eq: FastRisc}. To begin with, let $\alpha$ and $\varepsilon$ two small positive parameters. Our aim is first to obtain a suitable estimate for $\Delta\left( V\left({r}/{\alpha} \right) \right)$
in the region $\{ r \le \alpha \varepsilon \}$. We have:
\begin{equation}\label{eq: stimaLapU}
\Delta\left( V\left({r}/{\alpha} \right) \right) = -\frac{1}{\alpha^2(1-m)T}V^{1/m}(r/\alpha) + \frac{N-1}{\alpha^2}\left( \alpha \coth(r) - \coth(r/\alpha) \right)V^\prime(r/\alpha) \, .
\end{equation}
The function $h(r)=r \coth(r)$ is regular. In particular,
\begin{equation}\label{eq: stimaLapUTaylor}
h(r)= 1 + h^\prime(0) \, r + q(r) r^2  \, , \ \ h(r/\alpha)=1 + h^\prime(0) \frac{r}{\alpha} + q(r/\alpha) \left(\frac{r}{\alpha}\right)^2 \, ,
\end{equation}
where both $|q(r)|$ and $|q(r/\alpha)|$ can be bounded by
$$ Q= \max_{s \in [0,1]} \frac{h^{\prime\prime}(s)}{2} $$
provided $\alpha,\varepsilon$ are smaller than $1$. In order to control the right term on the r.h.s.\ of \eqref{eq: stimaLapU}, we use \eqref{eq: stimaLapUTaylor}:
\begin{equation*}
\begin{aligned}
\left| \frac{\alpha \coth(r) - \coth(r/\alpha) }{\alpha^2} \right| &= \left| \frac{r \coth(r) - \frac{r}{\alpha}\coth(r/\alpha) }{\alpha r} \right| \\ &= \left| \frac{h^\prime(0)}{\alpha} + q(r) \frac{r}{\alpha} - \frac{h^\prime(0)}{\alpha^2} - q(r/\alpha) \frac{r}{\alpha^3} \right|  \le \frac{C}{\alpha^2}
\end{aligned}
\end{equation*}
for a suitable constant $C>0$ independent of $\alpha,\varepsilon$. Notice that, since $V$ is regular and $V^\prime(0)=0$, there exists $D>0$ (independent of $ \alpha,\varepsilon $) such that $ \left|V^\prime(r/\alpha) \right| \le D \varepsilon $ for all $r \le \alpha\varepsilon $. Hence,
\begin{equation}\label{eq: stimaLapU2}
\left| \frac{N-1}{\alpha^2}\left( \alpha \coth(r) - \coth(r/\alpha) \right)V^\prime(r/\alpha) \right| \le \frac{(N-1)C D \varepsilon }{\alpha^2}  \, .
\end{equation}
Then recall that $ V(\cdot) $ is decreasing and $V(0)>0$, so that from \eqref{eq: stimaLapU} and \eqref{eq: stimaLapU2} we can claim that there exists a constant $E>0$ such that for any $\varepsilon>0$ sufficiently small (depending only on $V$, $m$ and $N$) there holds
\begin{equation}\label{eq: stimaLapUFinale}
\Delta\left( V\left({r}/{\alpha} \right) \right) \le -\frac{E}{\alpha^2T} \ \ \ \forall r \in (0,\alpha\varepsilon) \, .
\end{equation}
From the $ L^{q}$-$L^\infty $ smoothing effects (see \cite[Lem. 6.1]{DKV91} or \cite[Th. 4.1]{BGV}, together with \eqref{eq: stimeNorma0}) we know that there exist $A>0$ and $\tau_0>0$ such that $w(r,\tau) \le A$ for all $(r,\tau) \in  (0,\infty) \times (\tau_0,\infty)$. Let $A_1>A_2$ be two given positive constants and let $\varepsilon>0$ be so small that \eqref{eq: stimaLapUFinale} holds. For a fixed $\tau_\ast\ge \tau_0+1$ set $f(\tau)=(\tau-\tau_0)/(\tau^\ast-\tau_0)$. %In particular, $ f(t_0)=0 $, $ f(t^\ast)=1 $ and $ f^\prime(t) = 1/(t^\ast-t_0)$.
First we shall prove that if $\alpha>0$ is small enough then the function
\begin{equation*}\label{eq: barrieraZero}
\varphi(r,\tau)=\left[ A_1(1-f(\tau)) + A_2 f(\tau) \right]^{\frac{1}{m}} V(r/\alpha)^{\frac{1}{m}}
\end{equation*}
is a \emph{supersolution} to \eqref{eq: FastRisc} in the region
\begin{equation}\label{eq: barrieraZeroReg}
\left\{ (r,\tau) \in (0,\alpha\varepsilon) \times (\tau_0,\tau_\ast) \right\} .
\end{equation}
To this end, notice that
\begin{equation}\label{eq: barrieraZeroCond}
\varphi(r,\tau_0) = A_1^{\frac{1}{m}} V^{\frac{1}{m}}(r/\alpha) \ge A_1^{\frac{1}{m}} V^{\frac{1}{m}}(\varepsilon) \ \ \forall r \in (0,\alpha\varepsilon) \, , \ \  \varphi(\alpha\varepsilon,\tau) \ge A_2^{\frac{1}{m}} V^{\frac{1}{m}}(\varepsilon) \ \ \ \forall \tau \in (\tau_0,\tau_\ast) \, ,
\end{equation}
while derivatives of $\varphi$ give
\begin{equation} \label{eq: derTphi}
\begin{aligned}
\varphi_\tau(r,\tau) &= -\frac{1}{m} f^\prime(\tau)V^{\frac1m}(r/\alpha)(A_1-A_2) \left[ A_1(1-f(\tau)) + A_2 f(\tau)  \right]^{\frac{1}{m}-1} \\
 & \ge -\frac{1}{m} (A_1-A_2) A_1^{\frac{1}{m}-1} V^{\frac{1}{m}}(0)  \, ,
\end{aligned}
\end{equation}
\begin{equation*} \label{eq: laplHphi}
 \Delta(\varphi^m)(r,\tau) = \left[A_1(1-f(\tau)) + A_2 f(\tau)\right] \Delta(V(r/\alpha))  \le -\frac{A_2 E}{\alpha^2T}  \,
\end{equation*}
and
\begin{equation} \label{eq: stimaHphi}
\varphi(r,\tau) \le A_1^{\frac{1}{m}}V^{\frac{1}{m}}(0) \, .
\end{equation}
Collecting \eqref{eq: derTphi}-\eqref{eq: stimaHphi} we get that for $\varphi$ to be a supersolution in the region \eqref{eq: barrieraZeroReg} it is enough to ask
\begin{equation}\label{eq: condPhiFinale}
-\frac{1}{m} (A_1-A_2) A_1^{\frac{1}{m}-1} V^{\frac{1}{m}}(0) \ge -\frac{A_2 E}{\alpha^2T} + \frac{A_1^{\frac{1}{m}}V^{\frac{1}{m}}(0)}{(1-m)T} \, ,
\end{equation}
which is achieved by choosing $\alpha=\alpha(V,T,m,N,A_1,A_2)$ sufficiently small.

Now fix $\varepsilon>0$ sufficiently small. Set $A_1=A^m/V(\varepsilon)$, $A_2=V(0)/V(\varepsilon)$, where we assume without loss of generality that $A^m>V(0)$ and pick $\alpha(V,T,m,N,A,\varepsilon)$ that complies with \eqref{eq: condPhiFinale}. Thanks to \eqref{eq: barrieraZeroCond} we have
\begin{equation}\label{eq: barrieraZeroCondExpl}
\varphi(r,\tau_0)  \ge A \ \ \ \forall r \in [0,\alpha\varepsilon] \, , \ \  \varphi(\alpha\varepsilon,\tau) \ge V^{\frac{1}{m}}(0) \ \ \ \forall \tau \in [\tau_0,\tau_\ast] \, .
\end{equation}
%Notice that one is free to center $\varphi$ in any pole $x_0 \in \mathbb{H}^N$, provided one replaces $r$ with $\rho=d(x,x_0)$. Let us choose $x_0 \neq o $ so that the pole $o$ is contained in the open ball of radius $\alpha\varepsilon$ centered at $x_0$.
From the fact that $V(\cdot)$ is decreasing and from the local uniform convergence \eqref{eq: convLoc}, we can take $\tau_0$ so large that $w^m(\alpha\varepsilon,\tau) \le V(0)$ for all $\tau \ge \tau_0$ (notice that $\varphi$ is a supersolution independently of $\tau_0,\tau_\ast$ provided $\tau_\ast-\tau_0 \ge 1$). Since \eqref{eq: barrieraZeroCondExpl} holds we can conclude, by comparison, that $w \le \varphi$ in the region \eqref{eq: barrieraZeroReg}. In particular,
$$ w^m(r,\tau_\ast) \le \frac{V^{2}(0)}{V(\varepsilon)}  \ \ \ \forall r \in [0,\alpha\varepsilon] \, . $$
By the remarks above this last result is actually valid for \emph{all} $\tau_\ast \ge \tau_0+1$. Hence, since $ V(\varepsilon) \rightarrow V(0) > 0 $ as $ \varepsilon \rightarrow 0 $, we conclude that for any $\epsilon>0$ there exist $r_\epsilon$ so small and $\tau_\epsilon$ so large that \eqref{eq: boundIntornoZero} holds true.
\end{proof}
\end{lem}

Thanks to Lemmas \ref{lem: convLoc}, \ref{lem: posLocZero} and \ref{lem: boundLocZero} we can extend the result of Lemma \ref{lem: convLoc} down to ${r=0}$ and get the following
\begin{pro}\label{teo: convLocZero}
Let $w$ be the rescaled solution \eqref{eq: rescSol} to \eqref{eq: FastRisc} and $V$ the radial, positive energy solution to the stationary problem \eqref{eq: statProb}. Then
\begin{equation}\label{eq: convLocZero}
\lim_{\tau \to \infty}  \left\|  w^m(\tau)-V \right\|_{L^\infty_{\textnormal{loc}}([0,\infty))} = 0 \, ,
\end{equation}
that is $w^m(\tau)$ converges uniformly to $V$ in any compact set $ K \Subset [0,\infty)$ as $\tau \to \infty$.
\end{pro}
% GRILLO non vuole questo commento
%
%
%Starting from Lemma \ref{lem: convLoc} and from the uniform boundedness of $w$ (which holds thanks to smoothing effects as remarked above), the proof of Theorem \ref{teo: convLocZero} could have been simplified using the regularity results of \cite{DiBCont} (with minor modifications, they can be applied \emph{locally} to this context too, see also \cite{FS}), which give (locally) a uniform spatial modulus of continuity for $w(d(x,o),\tau)$. Hence the trajectory $\{ w(\tau) \}$ is precompact in $C_{loc}(\mathbb{H}^N)$.
%
%However, exploiting the fact that we are considering only radial solutions, we could avoid the use of such regularity results.

% cmq, per scongiurare dubbi: u_n converge ovunque a u, e non solo quasi, grazie alla positività locale, dove tutto è regolare; ad ogni modo in caso contrario si sistema tutto parlando di limiti essenziali
\section{Estimates from above} \label{Sec: Above}
The goal of this section is to bound the ratio $w^m(r,\tau)/V(r)$ in $ L^\infty((0,\infty)) $ (and not only in $ L^\infty_{\textnormal{loc}}([0,\infty)) $ as we did in Section \ref{Sec: asy}) from above. Since $V(r)$ behaves like $e^{-(N-1)\,r}$ at infinity (see \eqref{eq:asympt} or Lemma \ref{lem: SolStaz} below), it will be enough to give an upper bound for $w^m(r,\tau)/e^{-(N-1) \, r}$. In fact our goal is to prove the following result.

\begin{pro}\label{thm: stimaSopraFond}
Let $w$ be the rescaled solution associated to a nonnegative energy solution $u$ to \eqref{eq: FastRadial}, as in \eqref{eq: rescSol}. Then, for any $\varepsilon>0$, there exists a positive constant $Q^\prime = Q^\prime (u_0, m, N, \varepsilon)$ such that
\begin{equation}\label{eq: stimaSopraFond}
w(r,\tau) \le Q^\prime \, e^{-\frac{N-1}{m} \, r} \ \ \ \forall (r,\tau) \in [0,\infty) \times [\varepsilon,\infty) \, .
\end{equation}
\end{pro}

To this end, we begin with some preliminary lemmas.
\begin{lem}\label{lem: EmbLinf}
Let $v \in \Hr$. For any ${r_0}>0$ there exists a positive constant $C({r_0},N)$ such that
\begin{equation*}\label{eq: EmbLinf}
v(r) \leq  C({r_0},N) \, \| v \|_{H^1} \, e^{- \frac{N-1}{2} \, r }  \ \ \ \forall r \in [{r_0},\infty) \, .
\end{equation*}
\begin{proof}
Consider the function $z(r)=(\sinh r)^{N-1} \, v^2(r)$. We have:
$$ z^\prime(r)= (N-1) (\cosh r) (\sinh r)^{N-2} \, v^2(r) + 2 (\sinh r)^{N-1} \,  v(r) v^\prime(r) \, ; $$
integrating between $r$ and ${r_0}$ gives
\begin{equation}\label{eq: EmbLinf0}
\begin{split}
z(r) =& (\sinh {r_0})^{N-1} v^2({r_0}) + (N-1) \int_{{r_0}}^r{ v^2(s)   \, (\cosh s) (\sinh s)^{N-2} \, \mathrm{d}s } \\
     &+  2  \int_{{r_0}}^r{  v(s) v^\prime(s)   \, (\sinh s)^{N-1} \, \mathrm{d}s } \, .
\end{split}
\end{equation}
Since ${r_0}>0$ and the behaviour at infinity of $\sinh r$ and $\cosh r $ is the same, we can control the last two terms on the r.h.s. of \eqref{eq: EmbLinf0} with a constant (depending on ${r_0}$ and $N$) times $\| v \|_{H^1}^2$. As for the first term, notice that $H^1({r_0}/2,3r_0/2)$ is continuously embedded in $L^\infty({r_0}/2,3r_0/2)$ and $\Hr$ is in turn continuously embedded in $H^1({r_0}/2,3r_0/2)$ (again, through constants depending on ${r_0}$ and $N$). Hence, there exists $C({r_0},N)>0$ such that
$$ z(r) \leq C^2({r_0},N) \, \| v \|_{H^1}^2 \ \ \ \forall r \in [{r_0},\infty) \, , $$
which gives the claimed result since $\sinh r \asymp e^{r} $ for $r$ large.
\end{proof}
\end{lem}

\begin{lem} \label{lem: SolStaz}
For any $\widetilde{m} \in (m_s,1)$ there exists a solution $V$ to
\begin{equation}\label{eq: solStaz0}
-\Delta{V}(r)=V^{\frac{1}{\widetilde{m}}}(r) \ \ \ \forall r \in (0,\infty)
\end{equation}
which is smooth, strictly positive, belongs to $\Hr $ and satisfies
\begin{equation}\label{eq: solStazAnd}
A^{-1} \, e^{-(N-1)r} \le V(r) \le A \, e^{-(N-1)r} \ \ \ \forall r \in [0,\infty)
\end{equation}
for some positive constant $A=A(\widetilde{m})$.
\begin{proof}
As already mentioned in the Introduction, see \cite{MS08} and \cite{BGGV}.
\end{proof}
\end{lem}

\begin{lem} \label{lem: SopraSolBasic}
Let $u$ be a bounded, nonnegative energy solution to \eqref{eq: FastRadial}. There exists a positive constant $C_0=C_0(u_0,m,N)$ such that
\begin{equation}\label{eq: StimaSopraBasic}
u(r,t) \leq C_0 \, e^{-\frac{N-1}{2m} \, r} \ \ \ \forall(r,t) \in [0,\infty) \times (0,T) \, .
\end{equation}
\begin{proof}
It is a matter of straightforward computations to show that
\begin{equation}\label{eq: SopraBasic}
-\Delta\left( e^{-\frac{N-1}{2} \, r}  \right) \ge 0 \ \ \ \forall r \in (0,\infty) \, .
\end{equation}
Thanks to Lemma \ref{lem: EmbLinf}, the fact that $u_0^m \in \Hr$ and the boundedness of $u$, one can choose $C_0$ so large that $C_0 \, e^{-\frac{N-1}{2m} \, r}$ is above $u$ on a parabolic boundary of the type $[r_0,\infty) \times \{ 0 \} \cup \{ r_0 \} \times (0,T)$, for a given $r_0 \in (0, \infty)$. The conclusion then follows from \eqref{eq: SopraBasic}, the comparison principle and again the boundedness of $u$.
\end{proof}
\end{lem}
We are now ready to prove a better (spatial) estimate from above for nonnegative energy solutions to \eqref{eq: FastRadial}.
% Grillo: dire che vale anche per m sottocritico

\begin{lem}\label{lem: StimaSopra}
Let $u$ be a nonnegative energy solution to \eqref{eq: FastRadial}. For any ${t_\ast} \in (0,T)$ there exists a positive constant $Q=Q({t_\ast},u_0,m,N)$ such that
\begin{equation}\label{eq: StimaSopra}
u(r,{t_\ast}) \leq Q \, e^{-\frac{N-1}{m} \, r} \ \ \ \forall r \in [0,\infty) \, .
\end{equation}
\begin{proof}
We shall proceed by constructing a proper barrier. In particular, we shall prove that for a suitable choice of the parameter $\xi>0$, the following function is a supersolution to \eqref{eq: FastRadial} in the parabolic domain $(\xi, \infty) \times (0,t_\ast)$:
\begin{equation}\label{eq: costruzioneSopra}
\bar{u}(r,t)=C_0 \left[ A \, e^{\frac{N-1}{2} \xi}  \, V(r) \, f(t) + e^{-\frac{N-1}{2} r} \,  \left(1-f(t)\right)   \right]^{\frac{1}{m}} \, ,
\end{equation}
where $V$ is the solution to \eqref{eq: solStaz0} associated to a fixed $\widetilde{m} \in ({2m}/{(1+m)} , 1 )$ ($A$ being the corresponding constant that appears in \eqref{eq: solStazAnd}) and $f(t):[0,t_\ast] \rightarrow [0,1]$ is a regular increasing function such that $f(0)=0$, $f(t_\ast)=1$, which we shall define later. The constant $C_0$ is the one from \eqref{eq: StimaSopraBasic}: indeed, thanks to smoothing effects (recall the discussion in Section \ref{subs: pre}), there is no loss of generality in assuming that $u_0$ is bounded. Since
\begin{equation*} \label{eq: frontPar0}
\bar{u}(r,0) = C_0 \, e^{-\frac{N-1}{2m} \, r} \ge u_0(r) \ \ \ \forall r \in [0,\infty)
\end{equation*}
and
\begin{equation*} \label{eq: frontPar1}
\bar{u}(\xi,t) = C_0 \left[ A \, e^{\frac{N-1}{2} \xi}  \, V(\xi) \, f(t) + e^{-\frac{N-1}{2} \xi} \,  \left(1-f(t)\right)     \right]^{\frac{1}{m}} \ge  C_0 \, e^{-\frac{N-1}{2m} \, \xi} \ge  u(\xi,t)   \ \ \  \forall t \in [0,t_\ast) \, ,
\end{equation*}
in order to prove that $u(r,t) \leq \bar{u}(r,t) $ for all $(r,t) \in [\xi,\infty) \times [0,t_\ast]$ we are left with showing that, by choosing appropriately $\xi$, there holds
\begin{equation*}\label{eq: sopraEquazione}
\bar{u}_t(r,t) \geq \Delta\left( \bar{u}^m \right)(r,t)  \ \ \ \forall (r,t) \in (\xi,\infty) \times (0,t_\ast) \, .
\end{equation*}
We have:
\begin{equation*}\label{eq: derT}
\begin{split}
\bar{u}_t(r,t)&=C_0 \, f^\prime(t) \left[ A \, e^{\frac{N-1}{2} \xi}  \, V(r) - e^{-\frac{N-1}{2} r} \right] \frac{1}{m} \left[ A \, e^{\frac{N-1}{2} \xi}  \, V(r) \, f(t) + e^{-\frac{N-1}{2} r} \,  \left(1-f(t)\right)  \right]^{\frac{1}{m}-1}  \\
& \ge -C_0 \, \frac{1}{m} \, f^\prime(t) \, e^{-\frac{N-1}{2} r} \, \left[ A^2 \, e^{\frac{N-1}{2} \xi}  \, e^{-(N-1)r} \, f(t) + e^{-\frac{N-1}{2} r} \,  \left(1-f(t)\right)  \right]^{\frac{1}{m}-1}  \, .
\end{split}
\end{equation*}
Also, thanks to \eqref{eq: solStaz0}, % and \eqref{eq: solStazAnd},
\begin{equation*}\label{eq: laplacianR}
\begin{split}
\Delta\left( \bar{u}^m \right)(r,t) &= -C_0^m \left[ A \, e^{\frac{N-1}{2} \xi}  \, V^{\frac{1}{\widetilde{m}}}(r) \, f(t) + \frac{(N-1)^2}{4} \, e^{-\frac{N-1}{2} r} \, (2 \coth r -1) \, \left(1-f(t)  \right)    \right]  \\
&\le   -C_0^m \left[ A^{1-\frac{1}{\widetilde{m}}} \, e^{\frac{N-1}{2} \xi} \, e^{-\frac{N-1}{\widetilde{m}} r} f(t) + \frac{(N-1)^2}{4} \, e^{-\frac{N-1}{2} r} \, \left(1-f(t)  \right)  \right]  .
\end{split}
\end{equation*}
In particular, there exist two positive constants $B_0(m,\widetilde{m})$ and $B_1(\widetilde{m},N)$ such that
\begin{equation*}\label{eq: derT1}
-\bar{u}_t(r,t) \leq C_0 \, B_0  \,  f^\prime(t) \left[ e^{\frac{(N-1)(1-m)}{2m} \xi} \, e^{-\frac{(N-1)(2-m)}{2m} r}  f^{\frac{1}{m}-1}(t) +  e^{-\frac{N-1}{2m} r} \, \left(1-f(t)  \right)^{\frac{1}{m}-1}   \right]  ,
\end{equation*}
\begin{equation*}\label{eq: laplacianR0}
-\Delta\left( \bar{u}^m \right) (r,t) \geq C_0^m B_1  \left[ e^{\frac{N-1}{2} \xi} \, e^{-\frac{N-1}{\widetilde{m}} r} f(t) + e^{-\frac{N-1}{2} r} \, \left(1-f(t)  \right)  \right] .
\end{equation*}
Hence it is enough to show that, if $\xi$ is properly chosen, the following inequality holds in $(\xi,\infty) \times (0,t_\ast)$:
\begin{equation}\label{eq: diseqParab1}
\begin{split}
& C_0^{1-m} \, B_\ast \, f^\prime(t) \left[ e^{\frac{(N-1)(1-m)}{2m} \xi} \, e^{-\frac{(N-1)(2-m)}{2m} r}  f^{\frac{1}{m}-1}(t) +  e^{-\frac{N-1}{2m} r} \, \left(1-f(t)  \right)^{\frac{1}{m}-1}   \right]    \\
 \le & \, e^{\frac{N-1}{2} \xi} \, e^{-\frac{N-1}{\widetilde{m}} r} f(t) + e^{-\frac{N-1}{2} r} \, \left(1-f(t)  \right) ,
\end{split}
\end{equation}
where $B_\ast=B_\ast(m,\widetilde{m},N)$ is another suitable positive constant. Applying the change of variable $\rho=r-\xi$ and using the fact that $\widetilde{m}<1$ we infer that \eqref{eq: diseqParab1} is implied by
\begin{equation}\label{eq: diseqParab2}
\begin{split}
& C_0^{1-m} \, B_\ast \, f^\prime(t) \, e^{-\frac{N-1}{2m} \, \xi}  \left[ e^{-\frac{(N-1)(2-m)}{2m} \rho}  f^{\frac{1}{m}-1}(t) +  e^{-\frac{N-1}{2m} \rho} \, \left(1-f(t)  \right)^{\frac{1}{m}-1}   \right]   \\
\le & \, e^{-\frac{(N-1)(2-\widetilde{m})}{2\widetilde{m}} \xi} \left[ e^{-\frac{N-1}{\widetilde{m}} \rho} f(t) + e^{-\frac{N-1}{2} \rho}  \left( 1-f(t) \right)  \right] \ \ \ \forall (\rho,t) \in (0,\infty) \times (0,t_\ast)  \, .
\end{split}
\end{equation}
Since we choose $\widetilde{m}$ to lie in the interval $({2m}/{(1+m)},1 )$ we have that
$$\epsilon=(N-1)\left( \frac{1}{2m} - \frac{1}{\widetilde{m}} + \frac12 \right)>0 \, ; $$
therefore \eqref{eq: diseqParab2} reads
\begin{equation}\label{eq: diseqParab3}
C_0^{1-m} \, B_\ast \, e^{-\epsilon \, \xi} \, \underbrace{ f^\prime(t) \, \frac{e^{-\frac{(N-1)(2-m)}{2m} \rho}  f^{\frac{1}{m}-1}(t) +  e^{-\frac{N-1}{2m} \rho} \, \left(1-f(t)  \right)^{\frac{1}{m}-1} }   { e^{-\frac{N-1}{\widetilde{m}} \rho} f(t) + e^{-\frac{N-1}{2} \rho}  \left( 1-f(t) \right)  }  }_{l(\rho,t)} \le 1  \, .
\end{equation}
Our aim is now to show that for a suitable choice of $f(t)$ the function $l(\rho,t)$ remains bounded in $(0,\infty) \times (0,t_\ast)$. To this end, let us set $x=e^{-(N-1)\rho}$ and $f(t)=h(t/t_\ast)$, $h$ being a function to be defined which satisfies all the fulfilments required to $f$ in the interval $[0,1]$ instead of $ [0,t_\ast] $. The boundedness of $l(\rho,t)$ is implied by the boundedness of the ratio (recall that $\widetilde{m}>{2m}/{(1+m)}$)
\begin{equation}\label{eq: ratio1}
h^\prime(\tau) \frac{ x^{ \left( \frac{1}{m} - \frac12 \right) } h^{ \left( \frac{1}{m} -1 \right) }(\tau) + x^{\frac{1}{2m}} \left( 1-h(\tau) \right)^{\left(  \frac{1}{m}-1  \right) } }{ x^{\left(  \frac{1}{2m} + \frac12 \right) } h(\tau) + x^{\frac12} \left(1 - h(\tau) \right) } = h^\prime(\tau)  \frac{ x^\alpha h^\alpha(\tau) + x^{\frac{\alpha}{2}} \left( 1 -h(\tau) \right)^\alpha  }{ x^{\frac{\alpha+1}{2}}h(\tau)  + 1-h(\tau)   }  \, ,
\end{equation}
$$ \alpha=\frac{1}{m}-1 \, , $$
for $(x,\tau) \in (0,1) \times (0,1)$.  If $\alpha \ge 1 $, which corresponds to $m \le {1}/{2}$, the numerator in \eqref{eq: ratio1} is always smaller than or equal to the denominator. Therefore we remain with the case $\alpha \in (0,1)$, that is $m \in ({1}/{2},1) $. Here it is convenient to choose $h$ as
\begin{equation*}\label{eq: sceltaH}
h(\tau)=1-(1-\tau)^{\frac{1}{\alpha}} \, .
\end{equation*}
In this way, performing the change of variable $\sigma=1-\tau $, \eqref{eq: ratio1} becomes
\begin{equation}\label{eq: ratio2}
\frac{1}{\alpha} \, \sigma^{\frac{1}{\alpha}-1} \, \frac{x^\alpha\left( 1-\sigma^{\frac{1}{\alpha}} \right)^\alpha  + x^{\frac{\alpha}{2}} \sigma  }{x^{\frac{\alpha+1}{2}}\left(1-\sigma^{\frac{1}{\alpha}} \right) + \sigma^{\frac{1}{\alpha}}  }
\end{equation}
for $(x,\sigma) \in (0,1) \times (0,1)$. Now notice that, for $\sigma$ varying in $[1/2,1)$, \eqref{eq: ratio2} is bounded by a constant that depends only on $\alpha$ (recall that $\alpha \in (0,1)$). If instead $\sigma$ varies in $(0,1/2)$, the boundedness of \eqref{eq: ratio2} is equivalent to the boundedness of
\begin{equation}\label{eq: ratio3}
\frac{x^\alpha \sigma^{-1} + x^{\frac{\alpha}{2}}} {x^{\frac{\alpha+1}{2}} \sigma^{-\frac{1}{\alpha}}  +1} \le  \frac{x^\alpha \sigma^{-1} } {x^{\frac{\alpha+1}{2}} \sigma^{-\frac{1}{\alpha}}  +1}  + 1  \, .
\end{equation}
For any fixed $x \in (0,1)$, the maximum of the r.h.s. of \eqref{eq: ratio3} as $\sigma \in (0,\infty) $ can be found explicitly, and it is equal to
\begin{equation*}\label{eq: ratioMin}
x^{\frac{\alpha}{2}(1-\alpha)} {\alpha^\alpha}{\left(1-\alpha \right)^{1-\alpha}}  + 1 \, .
\end{equation*}
Summing up, we have proved that for a suitable choice of $f(t)$ (depending on whether $m \le 1/2$ or $m \in (1/2,1)$) the function $l(\rho,t)$ in \eqref{eq: diseqParab3} is bounded by a positive constant $K=K(t_\ast,m)$ as $(\rho,t)$ varies in $(0,\infty) \times (0,t_\ast)$. This means that if we set
$$ \xi \ge \frac{1}{\epsilon} \log{\left( C_0^{1-m} \, B_\ast  \, K \right)}  $$
we ensure that $\bar{u}(r,t)$ as in \eqref{eq: costruzioneSopra} is a supersolution to \eqref{eq: FastRadial} in the region $(\xi,\infty) \times (0,t_\ast) $. By the comparison principle, in particular,
\begin{equation*}\label{eq: stimaFinale}
u(r,t_\ast) \le \bar{u}(r,t_\ast) \le  C_0 \, A^{\frac{2}{m}} \, e^{\frac{N-1}{2m} \xi } \, e^{-\frac{N-1}{m} \, r}  \ \ \ \forall r \in (\xi,\infty) \, .
\end{equation*}
Since $u(\cdot,t_\ast)$ is also bounded in $(0,\xi]$, this gives \eqref{eq: StimaSopra}.
\end{proof}
\end{lem}
The result just proved is not enough in order to bound from above the ratio $w^m(r,\tau)/V(r)$ since it provides such boundedness only at any \emph{fixed} $\tau \in (0, \infty)$: indeed, recall that $u(\cdot, t)$ is bounded (so is $w(\cdot, \tau )$), $V(r)$ is locally bounded away from zero and its behaviour at infinity is the same as $e^{-(N-1) \, r}$. What estimate \eqref{eq: StimaSopra} lacks is a decay rate of order $(T - t_\ast)^{1/(1-m)}$ on the r.h.s.: we shall now prove that this is indeed the case.

\begin{lem}\label{lem: energyDec}
Let $w$ be the rescaled solution associated to a bounded, nonnegative energy solution $u$ to \eqref{eq: FastRadial}. There exists a constant $C_1 = C_1 (u_0,m,N)$ such that
\begin{equation}\label{eq: energyDec}
w(r, \tau ) \le C_1 \, e^{-\frac{N-1}{2m} \, r} \ \ \ \forall (r, \tau) \in [0, \infty) \times (0, \infty) \, .
\end{equation}
\begin{proof}
The fact that
\begin{equation}\label{eq: energyDec0}
\left\| u(t) \right\|_{\infty} \le D \left(T - t \right)^{\frac{1}{1-m}} \ \ \ \forall t \in (0,T)
\end{equation}
for a suitable positive constant $D=D(u_0,m,N)$ can be proved exactly as in \cite[Sects. 5-6]{DKV91} (see also the results of \cite{BGV}). Then notice that Lemma \ref{lem: stimeNorma} in particular yields
\begin{equation}\label{eq: energyDec1}
\left\| \left(w^m \right)^\prime \! (\tau) \right\|_{2} \le \left\| \left(u_0^m \right)^\prime \right\|_{2} \ \ \ \forall \tau \in (0,\infty) \, .
\end{equation}
The Poincar\'e inequality in \eqref{eq: poinSob}, Lemma \ref{lem: EmbLinf} and \eqref{eq: energyDec1} provide the following bound:
\begin{equation*}\label{eq: energyDec2}
w^m(r,\tau) \le D^\prime(r_0,N) \, \left\| \left(u_0^m \right)^\prime \right\|_{2} e^{-\frac{N-1}{2} \, r} \ \ \ \forall (r,\tau) \in [r_0,\infty) \times (0,\infty) \, ,
\end{equation*}
which together with \eqref{eq: energyDec0} gives the claimed estimate \eqref{eq: energyDec}.
\end{proof}
\end{lem}
Notice that we proved \eqref{eq: energyDec} under the hypothesis that $ u $ is a bounded, nonnegative energy solution. However, thanks to the aforementioned smoothing effects, it also holds for any solution associated to data as in the hypothesis of Theorem \ref{ThConvRel}, provided one starts from $ \tau=\varepsilon > 0$ instead of $\tau=0$.

The bound \eqref{eq: energyDec} for $w$ is the exactly the same as \eqref{eq: StimaSopraBasic} for $u$, which was a key starting point in order to prove the claim of Lemma \ref{lem: StimaSopra}. Indeed, as we shall see now, the barrier exploited in the proof of such Lemma also works for the equation solved by $w$. This allows us to obtain the next fundamental estimate from above.

\vskip8pt\noindent{\it Proof of Proposition \ref{thm: stimaSopraFond}}.
As just remarked, we need only prove that for a suitable choice of the parameter $\xi>0$ the barrier constructed in the proof of Lemma \ref{lem: StimaSopra} still works. So, for a given $\tau_\ast \in (\varepsilon,\infty)$, consider again the function
\begin{equation}\label{eq: costruzioneSopraFond}
\bar{w}(r,\tau)=C_1 \left[ A \, e^{\frac{N-1}{2} \xi}  \, V(r) \, f(\tau) + e^{-\frac{N-1}{2} r} \,  \left(1-f(\tau)\right)   \right]^{\frac{1}{m}}  ,
\end{equation}
where $\xi$, $V$, $A$, $f$ are as in \eqref{eq: costruzioneSopra} (just replace $t$ by $\tau$ and let $f(\varepsilon)=0$) and $C_1$ is the constant appearing in \eqref{eq: energyDec}. First of all, since $f(\cdot)$ is always included in $[0,1]$, $m < 1$ and $\widetilde{m} \in (2m/(1+m),1)$, we can bound the reaction term in \eqref{eq: FastRisc} (what makes \eqref{eq: FastRisc} actually differ from \eqref{eq: Fast}) in the following way:
\begin{equation*}\label{eq: stimaSopraFond0}
\begin{split}
\frac{1}{(1-m)T} \, \bar{w}(\rho+\xi,\tau) &\le C_1 \, B_2 \, e^{-\frac{N-1}{2m} \, \xi}  \left[ e^{-\frac{N-1}{m} \rho}  f^{\frac{1}{m}}(\tau) +  e^{-\frac{N-1}{2m} \rho} \, \left(1-f(\tau)  \right)^{\frac{1}{m}}   \right]  \\
& \le  C_1 \, B_2 \, e^{-\frac{N-1}{2m} \, \xi}  \left[ e^{-\frac{N-1}{\widetilde{m}} \rho}  f(\tau) +  e^{-\frac{N-1}{2} \rho} \, \left(1-f(\tau)  \right)  \right] \, ,
\end{split}
\end{equation*}
where $B_2$ is a positive constant depending only on $T$, $m$, $\widetilde{m}$ and we performed the usual change of space variable $\rho=r-\xi$. So, the equivalent of \eqref{eq: diseqParab2} in this context reads
\begin{equation*}\label{eq: stimaSopraFond1}
\begin{split}
& C_1 \, B_0 \, f^\prime(\tau) e^{-\epsilon \, \xi} \left[ e^{-\frac{(N-1)(2-m)}{2m} \rho}  f^{\frac{1}{m}-1}(\tau) +  e^{-\frac{N-1}{2m} \rho} \, \left(1-f(\tau)  \right)^{\frac{1}{m}-1}   \right]   \\
\le & \left( C_1^m B_1 - C_1 \, B_2 \,  e^{-\epsilon \, \xi}  \right) \left[ e^{-\frac{N-1}{\widetilde{m}} \rho}  f(\tau) +  e^{-\frac{N-1}{2} \rho} \, \left(1-f(\tau)  \right)  \right]
\end{split}
\end{equation*}
$$ \forall (\rho,\tau) \in (0,\infty) \times (\varepsilon,\tau_\ast) \, . $$
By elementary computations one gets that
$$ \frac{C_1^m B_1}{2} \le C_1^m B_1 - C_1 \, B_2 \,  e^{-\epsilon \, \xi}  $$
provided
\begin{equation}\label{eq: condExtraXi}
\xi \ge \frac{1}{\epsilon} \log{\left( \frac{2 C_1^{1-m} B_2}{B_1} \right)} \, .
\end{equation}
Therefore, under assumption \eqref{eq: condExtraXi}, we can repeat the same proof of Lemma \ref{lem: StimaSopra} starting from \eqref{eq: diseqParab3} (one replaces $B_\ast$ by $2 B_0 / B_1$ and $C_0$ by $C_1$). Hence we end up with the existence of a positive parameter $\xi=\xi(C_1(u_0,m,N,\varepsilon),T,\epsilon)$ such that
\begin{equation}\label{eq: stimaSopraFond2}
w(r,\tau) \leq \bar{w}(r,\tau) \ \ \ \forall (r,\tau) \in [\xi,\infty) \times [\varepsilon,\tau_\ast] \, .
\end{equation}
The validity of \eqref{eq: stimaSopraFond} is then a consequence of \eqref{eq: stimaSopraFond2}, \eqref{eq: costruzioneSopraFond} (evaluated at $\tau=\tau_\ast$) and \eqref{eq: energyDec0}. \qed

\begin{oss}\label{oss: StimaSopra} \rm
In the proofs of Lemmas \ref{lem: SopraSolBasic}, \ref{lem: StimaSopra} and Proposition \ref{thm: stimaSopraFond} we applied the comparison principle in parabolic regions of the form $(\xi,\infty)\times(0,t_\ast)$, without considering the $\{r=\infty \}$ side of the parabolic boundary. However, this technical issue is easy solvable by approximating $u$ with the solutions $u_n$ to \eqref{eq: Fastn}, applying comparison between $u_n$ and $\bar{u}$ in $(\xi,n) \times (0,t_\ast)$ and passing to the limit as $n \to \infty$ (using also the fact that $ T_n \uparrow T $).
\end{oss}
\begin{oss}\label{oss: AllM} \rm
The estimate from above provided by Lemma \ref{lem: StimaSopra} actually holds for all $m \in (0,1)$, since its method of proof only requires the existence of a solution to \eqref{eq: solStaz0} satisfying \eqref{eq: solStazAnd} for a value of $\widetilde m$ which can be taken as close to 1 as necessary.
\end{oss}

\section{Estimates from below}\label{sec: below}
This section is the dual of the previous one: our aim here is to bound the ratio $w^m(r,\tau)/V(r)$ from below. Again, thanks to \eqref{eq: solStazAnd}, this will be equivalent to providing a lower bound for $w^m(r,\tau)/e^{-(N-1)\,r}$. More precisely, we shall prove the following result.

\begin{pro}\label{thm: below}
Let $u$ be a nonnegative non-zero energy solution to \eqref{eq: FastRadial}. For any $\varepsilon \in (0,T)$ there exists a positive constant $P^\prime=P^\prime(u_0,m,N,\varepsilon)$ such that
\begin{equation}\label{eq: StimaSotto}
u(r,t_\ast) \ge P^\prime \, e^{-\frac{N-1}{m} \, r}  \left( T-t_\ast  \right)^{\frac{1}{1-m}}  \ \ \ \forall \left(r,t_\ast\right) \in [0,\infty) \times \left[\varepsilon , T \right)  .
\end{equation}
\end{pro}

To this end, we need a preliminary step.

\begin{lem}\label{lem: belowPre}
Let $u$ be a nonnegative non-zero energy solution to \eqref{eq: FastRadial}. For any given $\alpha > {N-1}$ and ${t_\ast} \in (0,T)$ there exists a positive constant $P=P({t_\ast},\alpha,u_0,m,N)$ such that
\begin{equation}\label{eq: StimaSottoPre}
u(r,{t_\ast}) \ge P \, e^{-\frac{\alpha}{m} \, r} \ \ \ \forall r \in [0,\infty) \, .
\end{equation}
\begin{proof}
We shall prove that for a suitable choice of the positive parameters $\mu_0$ and $\xi$ the following function is a subsolution to \eqref{eq: FastRadial} in the parabolic region $(\xi, \infty) \times (t_\ast/2,t_\ast) $:
\begin{equation}\label{eq: funcSottoPre}
\underline{u}(r,t) = \mu_0 \left[ \left(1 + e^{-\alpha \, (r-\xi)} \right) f(t) - 1  \right]_{+}^{\frac{1}{m}} \, ,
\end{equation}
where $f:[t_\ast/2,t_\ast] \rightarrow [1/2,1]$ is an increasing function such that $f(t_\ast/2)=1/2$, $f(t_\ast)=1$, to be defined later. We have:
\begin{equation*}\label{eq: StimaSottoPreDt}
\underline{u}_t(r,t) = \mu_0 \, \frac{1}{m} \, f^\prime(t) \left( 1 + e^{-\alpha \, (r-\xi)}  \right) \left[ \left(1 + e^{-\alpha \, (r-\xi)} \right) f(t) - 1  \right]^{\frac{1}{m}-1}
\end{equation*}
and
\begin{equation*}\label{eq: StimaSottoPreDxx}
\Delta \left( \underline{u}^m \right)(r,t) =  \mu_0^m \, f(t)  \, e^{-\alpha \, (r-\xi)}  \, \left(  \alpha^2  - (N-1) \, \alpha \coth r   \right)
\end{equation*}
in the region where $\left(1 + e^{-\alpha \, (r-\xi)} \right) f(t) - 1 $ is nonnegative (below, we shall always work tacitly in such region), while both ${\underline{u}}_t$ and $\Delta(\underline{u}^m)$ are zero outside it. Let us check conditions on the parabolic boundary. On $ [\xi,\infty) \times \left\{ t_\ast/2 \right\}$ $\underline{u}$ satisfies
\begin{equation*}\label{eq: SottoPrePB0}
\underline{u}(r,t_\ast/2) =  \mu_0 \left[ \frac12  \left(1 + e^{-\alpha \, (r-\xi)} \right) - 1  \right]_{+}^{\frac{1}{m}} = 0
\end{equation*}
and on  $\left\{ \xi \right\} \times \left(t_\ast/2, t_\ast \right) $ there holds
\begin{equation*}\label{eq: SottoPrePB1}
\underline{u}(\xi, t) =  \mu_0 \left[ 2 f(t) -1  \right]^{\frac{1}{m}} \le \mu_0 \, .
\end{equation*}
Therefore, in order to have $\underline{u} \le u $ on $[\xi,\infty) \times \left\{ t_\ast/2 \right\} \cup \left\{ \xi \right\} \times \left(t_\ast/2, t_\ast \right) $, we need only ask
\begin{equation*}\label{eq: SottoPrePB}
\mu_0 \le \inf_{t \in \left(\frac{t_\ast}{2},t_\ast \right)} u(\xi,t) =  \lambda(\xi,t_\ast,u_0)   > 0  \, ,
\end{equation*}
the last inequality following from standard positivity results (if $ t_\ast $ is close enough to $ T $, one can also exploit the results of Section \ref{Sec: asy}). Now let us check the differential equation. Upon the usual change of spatial variable $\rho=r-\xi $,
$$ \underline{u}_t(r,t) \le \Delta \left( \underline{u}^m \right)(r,t) \ \ \ \forall (r,t)  \in (\xi,\infty) \times \left( \frac{t_\ast}{2} , t_\ast \right) $$
reads
\begin{equation}\label{eq: SottoPreDis0}
\begin{split}
& \mu_0 \, \frac{1}{m} \, f^\prime(t) \, \left(1+e^{-\alpha \, \rho} \right) \left[ \left( 1+e^{-\alpha \, \rho} \right)f(t)  -1  \right]^{\frac{1}{m}-1}   \\
\le & \mu_0^m \, f(t)  \, e^{-\alpha \, \rho}  \, \left(  \alpha^2  - (N-1) \, \alpha \coth(\rho+\xi)   \right)   \ \ \ \forall (\rho,t) \in (0,\infty) \times \left( \frac{t_\ast}{2} , t_\ast \right) .
\end{split}
\end{equation}
If we choose $\xi=\xi(\alpha,N)$ so large that, for instance,
\begin{equation}\label{eq: SottoPreXi}
\alpha^2  - (N-1) \, \alpha \coth \xi \ge  \frac12 \left( \alpha^2  - (N-1) \, \alpha \right)=C_1(\alpha,N) > 0  \, ,
\end{equation}
we get that \eqref{eq: SottoPreDis0} is implied by
\begin{equation*}\label{eq: SottoPreDis1}
\mu_0^{1-m} \, \frac{1}{m \, C_1} \, f^\prime(t) \, \frac{  \left(1+e^{-\alpha \, \rho} \right) \left[ \left( 1+e^{-\alpha \, \rho} \right)f(t)  -1  \right]^{\frac{1}{m}-1}    }{ f(t)  \, e^{-\alpha \, \rho} }   \le 1  \ \ \ \forall (\rho,t) \in (0,\infty) \times \left( \frac{t_\ast}{2} , t_\ast \right) ,
\end{equation*}
which is in turn implied by (recall that $1/2 \leq f(t) \leq 1$)
\begin{equation}\label{eq: SottoPreDis2}
\mu_0^{1-m} \, \underbrace{ \frac{4}{m \, C_1} \, f^\prime(t) \, e^{\alpha \, \rho}  \, \left[ \left( 1+e^{-\alpha \, \rho} \right)f(t)  -1  \right]^{\frac{1}{m}-1}   }_{L(\rho,t)}    \le 1  \ \ \ \forall (\rho,t) \in (0,\infty) \times \left( \frac{t_\ast}{2} , t_\ast \right) .
\end{equation}
If $m \le 1/2$ the function $L$ in \eqref{eq: SottoPreDis2} is bounded from above (take $f$ regular) by a constant that depends only on $t_\ast$, $\alpha$, $m$ and $N$. If instead $m \in (1/2,1)$ this is in general false, unless one chooses $f$ carefully. To this end, consider the function
\begin{equation*}\label{eq: SottoPreDeff0}
h(\tau)=1-\frac12 \left(1 - \tau \right)^{\frac{m}{1-m}}    \ \ \ \forall \tau \in [0,1]
\end{equation*}
and set
\begin{equation*}\label{eq: SottoPreDeff1}
f(t) = h\left(2\frac{t}{t_\ast} - 1 \right)    \ \ \ \forall t \in \left[ \frac{t_\ast}{2} , t_\ast \right] .
\end{equation*}
Elementary computations (one can find the exact maximum of $L(\rho,t)$, in the region where $(1+e^{-\alpha \, \rho})f(t) \ge 1$, at any given $t$) show that
\begin{equation*}\label{eq: SottoPreL}
L(\rho,t) \le \frac{C_2}{t_\ast}  \ \ \  \forall (\rho,t) \in (0,\infty) \times \left(\frac{t_\ast}{2} , t_\ast \right)
\end{equation*}
for a suitable positive constant $C_2=C_2(\alpha,m,N)$ (which we assume to work for the case $m \le 1/2$ too). Hence, we proved that $\underline{u}$ as in \eqref{eq: funcSottoPre} is indeed a subsolution to \eqref{eq: FastRadial} providing that
\begin{equation*}\label{eq: SottoPreMu0}
\mu_0 \le \min\left( \lambda(\xi(\alpha,N),t_\ast,u_0) \, , \, \left(\frac{t_\ast}{C_2}\right)^{\frac{1}{1-m}} \right)  .
\end{equation*}
In particular, at $t=t_\ast$ there holds
$$ u(r,t_\ast) \ge \underline{u}(r,t_\ast) = \mu_0\left( {t_\ast}, \alpha, u_0, m, N \right) \, e^{\frac{\alpha}{m} \, \xi(\alpha,N)}  \, e^{-\frac{\alpha}{m} \, r }  \ \ \ \forall r \in (\xi,\infty)  \, , $$
which yields the thesis together with the local positivity of $u(\cdot,t_\ast)$ in $(0,\xi)$.
\end{proof}
\end{lem}
The result provided by the previous Lemma asserts that at any given $t_\ast \in (0,T(u_0))$ all nonnegative non-zero energy solutions of \eqref{eq: FastRadial} go to infinity (as $r \to \infty $) slower than $e^{-\frac{\alpha}{m} \, r}$, for any $\alpha > N-1$. Exploiting such fact, we are able to prove that one can actually take $\alpha=N-1$.

\vskip8pt\noindent{\it Proof of Proposition \ref{thm: below}}.
%For notational simplicity we shall prove the result in the time interval $((\varepsilon+T)/2,T)$. Indeed it is straightforward to modify the forthcoming proof so that it also applies to the time interval $(\varepsilon ,T)$, and this is enough to get the claim (up to admitting dependence of the multiplying constants on $ \varepsilon $ as well).
Once again we construct a lower barrier which has the desired property stated in \eqref{eq: StimaSotto}. Indeed, given $t_\ast \in (2\varepsilon ,T)$ (the final result will follow just by replacing $ \varepsilon $ with $\varepsilon/2 $), consider the function
\begin{equation}\label{eq: StimaSottoDefBar}
\underline{u}(r,t)= \mu_0 \left[ \left[\left( e^{-\beta \, (r-\xi)} + A^{-1} V(r) e^{(N-1) \, \xi} \right)f(t) - e^{-\beta \, (r-\xi)}  \right]_{+} + e^{-\alpha \, (r-\xi)} \right]^{\frac{1}{m}}  \, ,
\end{equation}
where $\alpha=\alpha(m,N)$ and $\beta=\beta(m,N)$ are fixed parameters such that
\begin{equation}\label{eq: StimaSottoCondAB}
\alpha > N-1 \, , \ \beta < N-1 \, , \ \alpha \le \beta + (N-1) \, \frac{1-m}{m} \, ,
\end{equation}
$f:[\varepsilon,t_\ast] \rightarrow [0,1]$ is a regular increasing function that satisfies $f(\varepsilon)=0$, $f(t_\ast)=1$ and $V$ is the solution to \eqref{eq: solStaz0} associated to $\widetilde{m}=m$ ($A$ being the corresponding constant appearing in \eqref{eq: solStazAnd}). We want to prove that, if one chooses the positive parameters $\xi$ and $\mu_0$ properly, then $\underline{u}$ is a subsolution to \eqref{eq: FastRadial} in the parabolic region $(\xi, \infty) \times (\varepsilon,t_\ast)$. By \eqref{eq: StimaSottoDefBar}, we have
\begin{equation*}\label{eq: StimaSottoPB0}
\underline{u}(r,t_0)= \mu_0 \, e^{\frac{\alpha}{m} \, \xi} \,  e^{-\frac{\alpha}{m} \, r}  \ \ \ \forall r \in [\xi,\infty)
\end{equation*}
and
\begin{equation*}\label{eq: StimaSottoPB1}
\underline{u}(\xi,t) \le  \mu_0 \left[ A^{-1} V(\xi) e^{(N-1) \, \xi} + 1 \right]^{\frac{1}{m}} \le \mu_0 \, 2^{\frac{1}{m}}  \ \ \ \forall t \in (\varepsilon,t_\ast)  \, .
\end{equation*}
Therefore $\underline{u}$ and $u$ are correctly ordered on the parabolic boundary $[\xi,\infty) \times \left\{ \varepsilon \right\} \cup \left\{ \xi \right\} \times (\varepsilon,t_\ast) $ provided (recall \eqref{eq: StimaSottoPre})
$$ \mu_0 \le e^{-\frac{\alpha}{m} \, \xi} \, P(\varepsilon,\alpha,u_0,m,N) \textnormal{ \ \ \ and \ \ \ }  \mu_0 \le 2^{-\frac{1}{m}} \, \inf_{t \in (\varepsilon,t_\ast)} u(\xi,t) = \lambda(\xi,t_\ast,u_0,m,N,\varepsilon) >0 \, , $$
that is
\begin{equation}\label{eq: StimaSottoMu0}
\mu_0 \le \min \left( e^{-\frac{\alpha}{m} \, \xi} \, P \, , \, \lambda \right) .
\end{equation}
Now let us compute the derivatives of $\underline{u}(r,t)$:
\begin{equation*}\label{eq: StimaSottoDt}
\begin{split}
\underline{u}_t (r,t) = &  \mu_0 \, \frac{1}{m} \, f^\prime(t) \left( e^{-\beta \, (r-\xi)} + A^{-1} V(r) e^{(N-1) \, \xi} \right) \, \sign_{+} q(r,t)  \\
& \times \left[ \left[\left( e^{-\beta \, (r-\xi)} + A^{-1} V(r) e^{(N-1) \, \xi} \right)f(t) - e^{-\beta \, (r-\xi)}  \right]_{+} + e^{-\alpha \, (r-\xi)} \right]^{\frac{1}{m}-1}  ,
\end{split}
\end{equation*}
\begin{equation*}\label{eq: StimaSottoDxx}
\begin{split}
  & \Delta\left( \underline{u}^m \right)(r,t)  \\
 = & \mu_0^m  \left[ \left[ \left( f(t)-1 \right) (\beta^2-(N-1) \, \beta \coth r ) \, e^{-\beta \, (r-\xi)} - A^{-1} V^{\frac{1}{m}}(r) e^{(N-1) \, \xi} f(t) \right] \, \sign_{+} q(r,t)  \right. \\
& \phantom{=\mu_0^m  \left[ \right. } + \left. (\alpha^2-(N-1) \, \alpha \coth r ) \, e^{-\alpha \, (r-\xi)} \right] ,
\end{split}
\end{equation*}
where for the sake of notational convenience we have set
$$ q(r,t) = \left( e^{-\beta \, (r-\xi)} + A^{-1} V(r) e^{(N-1) \, \xi} \right)f(t) - e^{-\beta \, (r-\xi)} \, . $$
If we take $\xi=\xi(\alpha,N)$ large enough so that \eqref{eq: SottoPreXi} holds, set $\rho=r-\xi$ and use \eqref{eq: solStazAnd}, we obtain:
\begin{equation*}\label{eq: StimaSottoDt1}
\underline{u}_t (\rho+\xi,t) \le  \mu_0 \, \frac{1}{m} \, f^\prime(t) \left( e^{-\beta \, \rho} +  e^{-(N-1) \, \rho} \right) \! \left[ \left[\left( e^{-\beta \, \rho} + e^{-(N-1) \, \rho} \right)f(t) - e^{-\beta \, \rho}  \right]_{+} + e^{-\alpha \, \rho} \right]^{\frac{1}{m}-1} \! ,
\end{equation*}
\begin{equation*}\label{eq: StimaSottoDxx1}
\Delta\left( \underline{u}^m \right)(\rho+\xi,t) \ge  \mu_0^m  \left[-A^{\frac{1-m}{m}} \, e^{-(N-1)\frac{1-m}{m} \xi} \, f(t) \, e^{-\frac{N-1}{m} \, \rho}   + C_1 \, e^{-\alpha \, \rho} \right]
\end{equation*}
\begin{equation*}
\forall (\rho,t) \in (0,\infty) \times (\varepsilon,t_\ast) \, .
\end{equation*}
Hence once we choose $\xi=\xi(\alpha,m,N)$ such that, in addition to \eqref{eq: SottoPreXi}, it also satisfies
$$ A^{\frac{1-m}{m}} \, e^{-(N-1)\frac{1-m}{m} \xi} \le \frac{C_1}{2} \, , $$
we deduce that in order to have $\underline{u}\le u$ in $(\xi,\infty)\times(\varepsilon,t_\ast)$ it is enough to ask (recall that $f(t) \in [0,1]$)
\begin{equation}\label{eq: StimaSottoDisFin0}
\mu_0 \, \frac{1}{m} \, f^\prime(t) \left( e^{-\beta \, \rho} +  e^{-(N-1) \, \rho} \right) \left[ e^{-(N-1) \, \rho} + e^{-\alpha \, \rho} \right]^{\frac{1}{m}-1} \le \mu_0^m \, \frac{C_1}{2} \, e^{-\alpha \, \rho}
\end{equation}
$$ \forall (\rho,t) \in (0,\infty) \times (\varepsilon,t_\ast) \, . $$
Upon setting $f(t)=h((t-\varepsilon)/(t_\ast-\varepsilon))$, $h:[0,1]\rightarrow[0,1]$ being a given regular, increasing function such that $h(0)=0$ and $h(1)=1$, \eqref{eq: StimaSottoDisFin0} reads
\begin{equation}\label{eq: StimaSottoDisFin1}
\mu_0^{1-m} \, \underbrace{ \frac{2}{m \, C_1 \, (t_\ast-\varepsilon)} \, h^\prime\left(\frac{t-\varepsilon}{t_\ast-\varepsilon}\right) \left( e^{(\alpha-\beta) \, \rho} +  e^{(\alpha-N+1) \, \rho} \right) \left[ e^{-(N-1) \, \rho} + e^{-\alpha \, \rho} \right]^{\frac{1}{m}-1} }_{L(\rho,t)} \le 1
\end{equation}
$$ \forall (\rho,t) \in (0,\infty) \times (\varepsilon,t_\ast) \, . $$
Thanks to \eqref{eq: StimaSottoCondAB} and to the fact that $t_\ast \in (2\varepsilon,T)$, the function $L(\rho,t)$ in \eqref{eq: StimaSottoDisFin1} is bounded in $(0,\infty)\times(\varepsilon,t_\ast)$ by a positive constant $C_2=C_2(T,m,N,\varepsilon)$. Therefore $\underline{u}(r,t)$ is a subsolution to \eqref{eq: FastRadial} providing that (recall \eqref{eq: StimaSottoMu0})
\begin{equation*}\label{eq: StimaSottoFinMu0}
\mu_0  \le \min\left( e^{-\frac{\alpha}{m} \, \xi} \, P \ , \ C_2^{-\frac{1}{1-m}} \ , \  2^{-\frac{1}{m}} \, \inf_{t \in (\varepsilon,t_\ast)} u(\xi,t) \right) ,
\end{equation*}
which is implied by
\begin{equation}\label{eq: StimaSottoFinMu1}
\mu_0  \le C_3 \min\left(1  \, , \, \inf_{t \in (\varepsilon,t_\ast)} u(\xi,t) \right)
\end{equation}
% Grillo: dire che fino a qui sostanzialmente ho una stima dal basso, fissato t, anche per m sottocritico (anche come remark alla fine di entrabe le stime, alto e basso)
for a suitable positive constant $C_3=C_3(\xi,u_0,m,N,\varepsilon)$. Exploiting Proposition \ref{teo: convLocZero} we can give a quantitative lower bound for the r.h.s. of \eqref{eq: StimaSottoFinMu1}. Indeed \eqref{eq: convLocZero} yields (in particular) the existence, for any given $\xi>0$, of a time $\hat{t}=\hat{t}(\xi,u_0)$ such that
\begin{equation}\label{eq: StimaSottoFinMu2}
u(\xi,t) \ge \frac{V^{\frac{1}{m}}(\xi)}{2} \, \left(T-t\right)^{\frac{1}{1-m}} \ \ \ \forall t \in \left( \hat{t} , T \right) .
\end{equation}
From standard positivity results we also know that, should $\varepsilon < \hat{t}$, $u(\xi,t)$ is still positive between $\varepsilon$ and $\hat{t}$; this fact and \eqref{eq: StimaSottoFinMu2} (together with the local positivity of $V$) ensure the existence of a positive constant $C_4=C_4(\xi,u_0,m,N,\varepsilon)$ such that
\begin{equation*}\label{eq: StimaSottoFinMu3}
u(\xi,t) \ge C_4 \, \left(T-t\right)^{\frac{1}{1-m}} \ \ \ \forall t \in \left(\varepsilon , T \right) ,
\end{equation*}
which gives
\begin{equation}\label{eq: StimaSottoFinMu4}
\inf_{t \in (\varepsilon,t_\ast)} u(\xi,t) \ge C_4 \, \left(T-t_\ast \right)^{\frac{1}{1-m}} \, .
\end{equation}
Combining \eqref{eq: StimaSottoFinMu4} and \eqref{eq: StimaSottoFinMu1} we infer that \eqref{eq: StimaSottoFinMu1} is implied by
\begin{equation}\label{eq: StimaSottoFinMu5}
 \mu_0 \le C_5 \, \left(T-t_\ast \right)^{\frac{1}{1-m}}
\end{equation}
for another positive constant $C_5=C_5(\xi,u_0,m,N,\varepsilon)$. The validity of \eqref{eq: StimaSotto} for $r$ varying in $(\xi,\infty)$ is then a consequence of choosing $\mu_0$ as the r.h.s. of \eqref{eq: StimaSottoFinMu5} and evaluating the subsolution $\underline{u}(r,t)$ at $t=t_\ast$, recalling \eqref{eq: solStazAnd}. On the contrary, the validity of \eqref{eq: StimaSotto} as $r$ varies in $(0,\xi]$ is a direct consequence of \eqref{eq: convLocZero} and the local positivity of $u$. \qed
\vskip10pt
\noindent Proposition \ref{thm: below} can of course be reformulated as follows.
\begin{cor}\label{cor: below}
Let $w$ be the rescaled solution associated to a positive energy solution $u$ to \eqref{eq: FastRadial}. For any ${\tau_0}>0$ there exists a positive constant $P^{\prime\prime} = ({\tau_0},u_0, m, N)$ such that
\begin{equation*}\label{eq: stimaSottoFond}
w(r,\tau) \ge P^{\prime\prime} \, e^{-\frac{N-1}{m} \, r} \ \ \ \forall (r,\tau) \in [0,\infty) \times [{\tau_0},\infty) \, .
\end{equation*}
\end{cor}

Notice that, as already anticipated in the Introduction, from Propositions \ref{thm: stimaSopraFond} and \ref{thm: below} we deduce formula \eqref{eq: ghp} in Theorem \ref{CorHarnack}.

\begin{oss}\label{oss: epsilon1}\rm
Both in the proofs of Lemma \ref{lem: belowPre} and Proposition \ref{thm: below}, when computing the hyperbolic Laplacian of the barriers \eqref{eq: funcSottoPre} and \eqref{eq: StimaSottoDefBar}, we neglected Dirac terms coming out from the second derivatives of positive parts. However this makes no problem since it is easy to check that such terms are nonnegative.
\end{oss}
\begin{oss}\label{oss: epsilon2}\rm
When we applied the comparison principle to $\underline{u}$ and $u$, both in the proofs of Lemma \ref{lem: belowPre} and Proposition \ref{thm: below}, we did not take into account $\{r=\infty\}$ as part of the parabolic boundary. To justify more rigorously those passages, it is enough to consider the following family of modified barriers:
\begin{equation*}\label{eq: modBar}
\underline{u}_\epsilon= \left[ \underline{u}^m - \epsilon \right]_{+}^{\frac{1}{m}} \le \underline{u} \, , \ \ \epsilon>0 \, ,
\end{equation*}
where $\underline{u}$ is either \eqref{eq: funcSottoPre} or \eqref{eq: StimaSottoDefBar}. Straightforward computations show that $\underline{u}_\epsilon$ is a subsolution to \eqref{eq: FastRadial} as long as $\underline{u}$ is. Moreover, for any fixed $\epsilon>0$, $\underline{u}_\epsilon(\cdot,t)$ has the property of being zero outside a compact set of the form $[\xi,R(\epsilon)] \subset [\xi,\infty)$. One then applies the comparison principle in $(\xi,R(\epsilon)) \times (\varepsilon,t_\ast) $ (let $\varepsilon=t_\ast/2$ when $\underline{u}$ is as in \eqref{eq: funcSottoPre}), gets $\underline{u}_\epsilon \le u $ and lets $\epsilon \rightarrow 0$.
\end{oss}
\begin{oss}\label{oss: AllMBelow} \rm
The estimate from below \eqref{eq: StimaSotto} holds for all $m \in (0,1)$, even though one has to admit that the constant $P^\prime$ in there, for $m \in (0,m_s]$, depends on $t^\ast$ as well. Indeed, in the proof of Proposition \ref{thm: below}, for simplicity we exploited the existence of a solution $V$ to \eqref{eq: solStaz0} satisfying \eqref{eq: solStazAnd} for $\widetilde{m}=m$. However, if $m \in (0,m_s]$, such a solution does not exist. Nonetheless it is easy to see that the choice of any $\widetilde{m} \in (m_s,1)$ instead of $m$ would have worked the same (provided one requires in addition that $\alpha < (N-1)/\widetilde{m}$).

Hence, this result together with the one discussed in Remark \ref{oss: AllM} prove the bounds \eqref{eq: sub} of Theorem \ref{CorSub}.
\end{oss}

\section{Results for the relative error}\label{Sec: rel}
As previously discussed, our interest is to show that the solution $u$ to \eqref{eq: FastRadial} converges in a strong sense to the stationary state $V$ which solves \eqref{eq: statProb}. More precisely, in this Section we shall prove \emph{uniform convergence in relative error}, namely
\begin{equation}\label{eq: limErrRel}
\lim_{\tau \rightarrow \infty} \left\| \frac{w^m(\tau)}{V}-1 \right\|_\infty = 0 \, ,
\end{equation}
that is \eqref{eq: ConvRelFin} rewritten in rescaled variables.

Since $V$ is strictly positive in any compact subset of $\mathbb{H}^N$, as a direct consequence of Proposition \ref{teo: convLocZero} (formula \eqref{eq: convLocZero}) we already know that \eqref{eq: limErrRel} holds \emph{locally} on $\mathbb{H}^N$. The nontrivial point is to prove that \eqref{eq: convLocZero} holds up to $r=\infty$, what we are concerned with in this section.

To begin with, let us write the equation solved by the \emph{relative error} $\phi=w^m/V - 1$. %Exploiting the property
%$$ \Delta\left(f \cdot g\right) = g \cdot \Delta(f) + f \cdot \Delta(g) + 2 f^\prime g^\prime $$
%and
Using the fact that $V$ satisfies \eqref{eq: statProb}, we get:
\begin{equation}\label{eq: PiuErrRel}
\frac{1}{m} \left(1+\phi \right)^{\frac{1}{m}-1} \phi_\tau = V^{1-\frac{1}{m}}\Delta \phi + 2 \, \frac{V^\prime \phi^\prime}{V^{\frac{1}{m}}} + \frac{1}{(1-m)T} \left[ \left(1+\phi \right)^{\frac{1}{m}} - \left(1+\phi \right) \right] ,  % \ \ \ \forall(r,t) \in (0,\infty) \times (1,\infty)
\end{equation}
where, as remarked above, the apex $^\prime$ stands for derivation w.r.t.\ $r$. It is also important to consider the equation solved by \emph{minus} the relative error $\psi=1-w^m/V$, which is
\begin{equation*}\label{eq: MenErrRel}
\frac{1}{m} \left(1-\psi \right)^{\frac{1}{m}-1} \psi_\tau = V^{1-\frac{1}{m}}\Delta \psi + 2 \, \frac{V^\prime \psi^\prime}{V^{\frac{1}{m}}} + \frac{1}{(1-m)T} \left[ \left(1-\psi \right) - \left(1-\psi \right)^{\frac{1}{m}} \right] .  % \ \ \ \forall(r,t) \in (0,\infty) \times (1,\infty)
\end{equation*}
In order to prove \eqref{eq: limErrRel} we shall construct suitable upper barriers both for $\phi$ and $\psi$, following the approach of \cite{BGVFast}.

The next lemma shows a good approximation property for $V$, which will turn out to be very useful to overcome some technical difficulties related to the upper barrier for $\psi$.
\begin{lem}\label{lem: Approx}
There exists a sequence of positive, regular radial solutions $\{V_n\}$ to
% i minimi sobolev su B_n dovrebbero essere unici e radiali, comunque eventualmente uno può risolvere i minimi col vincolo di essere radiali, e tutto funziona
\begin{equation}\label{eq: statProbN}
\begin{cases}
-\Delta V_n = \frac{1}{(1-m)T_n} \, V_n^{\frac{1}{m}}   & \textnormal{on ${B}_{n}$}  \\
  V_n=0    & \textnormal{on $\partial{B}_{n}$}
\end{cases}  \, ,
\end{equation}
where ${B}_{n}$ is the hyperbolic ball of radius $n$ centered at $o$, such that:
\begin{itemize}
\item $V_n(r) \le V(r)$ for all $r \ge 1$ and $ V_n \rightarrow V$ pointwise;
\item $ T_n \rightarrow T$;
\item for $\epsilon>0$ arbitrarily small, one can choose $n_{\epsilon}$ and $r_\epsilon$ so large that the following inequality holds:
\begin{equation}\label{eq: statProbNDeriv}
V_n^\prime(r) \le (1-\epsilon) \, V^\prime(r) \ \ \ \forall r \in (r_\epsilon,n) \, , \ \forall n \ge n_\epsilon \, .
\end{equation}
\end{itemize}
\begin{proof}
For any given $n \in \mathbb{N}$, let $X_n$ be the set of all functions $v \in \Hr$ such that $v(r)=0$ for all $r \ge n$ and $\| v \|_{\frac{m+1}{m}}=\| V \|_{\frac{m+1}{m}}$. Consider the following minimization problem:
\begin{equation}\label{eq: probMin}
\textnormal{Find $v_n \in X_n$:} \ \ \ \left\| v_n^\prime \right\|_{2} = \min_{v \in X_n} \left\| v^\prime \right\|_{2} \, .
\end{equation}
Thanks to standard arguments (in particular, the compact embedding of $\Hr$ into $L^{(m+1)/m}_{\textnormal{rad}}$), the solution $v_n$ of \eqref{eq: probMin} exists, is unique and solves \eqref{eq: statProbN} for some $T_n>0$. Since $ X_n \subset X_{n+1}$, the numerical sequence $\{ {\left\| v^\prime_n \right\|_{2}} \}$ is non-increasing. This implies that $\{ T_n\}$ is non-decreasing: indeed multiplying \eqref{eq: statProbN} by $v_n$ (upon replacing $V_n$ by $v_n$) and integrating by parts yield
$$ T_n =  \frac{\left\| v_n \right\|_{\frac{m+1}{m}}^{\frac{m+1}{m}}}{(1-m) \left\| v_n^\prime \right\|_{2}^2} = \frac{\left\| V \right\|_{\frac{m+1}{m}}^{\frac{m+1}{m}}}{(1-m) \left\| v_n^\prime \right\|_{2}^2}  \, . $$
Let us call $\widetilde{T}$ the limiting value of $\{ T_n\}$. Because $\{ {\left\| v^\prime_n \right\|_{2}} \}$ is bounded, $\{v_n\}$ (along subsequences) converges weakly in $\Hr$ and therefore strongly in $L^{{(m+1)}/{m}}_{\textnormal{rad}}$ to a certain function $\widetilde{V}$. Passing to the limit in the (weak formulation of the) equation solved by $v_n$, we get:
\begin{equation}\label{eq: statProbNLimit}
-\Delta \widetilde{V} = \frac{1}{(1-m)\widetilde{T}} \, \widetilde{V}^{\frac{1}{m}}   \ \ \ \textnormal{on $\mathbb{H}^N$} \, .
\end{equation}
First of all note that $\widetilde{T}$ cannot be $\infty$: if it were, from the Poincar\'e inequality in \eqref{eq: poinSob}  $\widetilde{V}$ would be zero, while we know that
$$\| \widetilde{V} \|_{{(m+1)}/{m}}= \lim_{n \rightarrow \infty} \| v_n \|_{{(m+1)}/{m}} =\| V \|_{{(m+1)}/{m}} > 0  \, .$$
So $\widetilde{V}$ is a positive, energy solution to \eqref{eq: statProbNLimit} for some $\widetilde{T} < \infty$, having the same $L^{{(m+1)}/{m}}_{\textnormal{rad}}$ norm as $V$: by the uniqueness result recalled in the beginning of Section \ref{Sec: asy}, it necessarily coincides with $V$ (and so $T_n \uparrow T$).

The next step is to prove that $\{v_n\}$ converges to $V$ also in $L^{{1}/{m}}_\textnormal{rad}$. Since $\{ v_n \}$ converges pointwise to $V$ (up to subsequences), and locally $L^{{(m+1)}/{m}}_\textnormal{rad}$ is continuously embedded in $L^{{1}/{m}}_\textnormal{rad}$, we need only dominate $\{ v_n\}$ outside a compact set with a function that belongs to $L^{{1}/{m}}_\textnormal{rad}$. To this end note that, from Lemma \ref{lem: EmbLinf} and from the fact that $\{ \left\| v^\prime_n \right\|_{2} \}$ is non-increasing, there exists a constant $K>0$ such that
\begin{equation}\label{eq: EmbLinfLemma}
v_n(r) \leq  K \, e^{- \frac{N-1}{2} \, r }  \ \ \ \forall r \in [{1},\infty) \, , \ \forall n \, .
\end{equation}
When $m$ is greater than or equal to $1/2$, the function on the r.h.s.\ of \eqref{eq: EmbLinfLemma} does not belong to $L^{{1}/{m}}_\textnormal{rad}$. Therefore we have to improve \eqref{eq: EmbLinfLemma} using the equation solved by $v_n$: integrating it between $0$ and $r$ we obtain the equality
\begin{equation}\label{eq: fundEqInt}
-(\sinh r)^{N-1} \, v_n^\prime(r) = {\int}_0^r \frac{v_n^{\frac{1}{m}}(s)}{T_n} \, (\sinh s)^{N-1} \mathrm{d}s \, .
\end{equation}
Suppose now that \eqref{eq: EmbLinfLemma} holds with a generic exponent $-(N-1)a$ (let $a>0$, $a \neq m$) in place of $-(N-1)/2$. Exploiting the corresponding analogues of \eqref{eq: EmbLinfLemma} and \eqref{eq: fundEqInt}, after some straightforward computations we get:
\begin{equation}\label{eq: EmbLinfLemma1}
-v_n^\prime (r) \leq  Q \, \frac{m}{m-a}  \, e^{-(N-1) \, r } \left(e^{(N-1)\frac{m-a}{m} \, r} -1\right) \ \ \ \forall r \in [{1},\infty) \, ,
\end{equation}
where $Q>0$ is a suitable positive constant that does not depend on $n$ (which may change from line to line). Since $v_n(n)=0$, we have:
\begin{equation}\label{eq: EmbLinfLemma2}
v_n(r) = \int_r^{n} -v_n^\prime(s) \, \mathrm{d}s \ \ \ \forall r \in [{1},n) \, .
\end{equation}
If $a > m$  \eqref{eq: EmbLinfLemma1} gives $-v_n^\prime (r) \leq  Q \, e^{-(N-1) \, r }$, in which case an integration of \eqref{eq: EmbLinfLemma2} entails
\begin{equation}\label{eq: EmbLinfLemma3}
v_n (r) \leq  Q \, e^{-(N-1) \, r } \ \ \ \forall r \in [{1},\infty) \, .
\end{equation}
If instead $a < m$ \eqref{eq: EmbLinfLemma1} gives $-v_n^\prime (r) \leq  Q \, e^{-(N-1) \, \frac{a}{m} \, r }$, and still integrating \eqref{eq: EmbLinfLemma2} one gets
\begin{equation*}\label{eq: EmbLinfLemma4}
v_n (r) \leq  Q \, e^{-(N-1) \, \frac{a}{m} \, r } \ \ \ \forall r \in [{1},\infty)
\end{equation*}
(the case $a=m$ can be dealt with similarly). It is plain that starting from $a=1/2$ and proceeding in this way, after a finite number of steps we obtain \eqref{eq: EmbLinfLemma3}. Since $e^{-(N-1) \, r }$ belongs to $L^{{1}/{m}}_\textnormal{rad}$, we have our dominating function and the convergence of $\{v_n\}$ to $V$ takes place also in $L^{{1}/{m}}_\textnormal{rad}$. Such convergence is crucial because it ensures that \eqref{eq: statProbNDeriv} holds for the sequence $\{ v_n \}$. Indeed applying \eqref{eq: fundEqInt} to $v_n=V$ yields
\begin{equation}\label{eq: fundEqIntLim}
\lim_{r \rightarrow \infty} -(\sinh r)^{N-1} \, V^\prime(r) =  \int_0^\infty \, \frac{V^{\frac{1}{m}}(s)}{T} \, (\sinh s)^{N-1} \mathrm{d}s  = C_V >0 \, ;
\end{equation}
for a given $\epsilon>0$, take $r_\epsilon$ so large that
\begin{equation*}\label{eq: fundEqIntLimEps}
\int_{r_\epsilon}^\infty \, \frac{V^{\frac{1}{m}}(s)}{T} \, (\sinh s)^{N-1} \mathrm{d}s  \le C_V\frac{\epsilon}{2} \, .
\end{equation*}
Thanks to the just proved convergence of $\{v_n\}$ to $V$ in $L^{{1}/{m}}_\textnormal{rad}$ and to the fact that $ T_n \uparrow T$, there exists $n=n(r_\epsilon,\epsilon)$, sufficiently large, such that
\begin{equation*}\label{eq: fundEqIntLimEps0}
\begin{aligned}
-(\sinh r)^{N-1} \, v_n^\prime(r) & = {\int}_0^r \frac{v_n^{\frac{1}{m}}(s)}{T_n} \, (\sinh s)^{N-1} \mathrm{d}s \ge {\int}_0^{r_\epsilon} \frac{v_n^{\frac{1}{m}}(s)}{T_n} \, (\sinh s)^{N-1} \mathrm{d}s \\
& \ge {\int}_0^{r_\epsilon} \frac{V^{\frac{1}{m}}(s)}{T} \, (\sinh s)^{N-1} \mathrm{d}s - C_V\frac{\epsilon}{2} \ge C_V (1-\epsilon)
\end{aligned}
\end{equation*}
$$ \forall r \in (r_\epsilon,n) \, , \ \forall n \ge n(r_\epsilon,\epsilon) \, . $$
Hence there holds
\begin{equation*}\label{eq: fundEqIntLimEps1}
\frac{-v_n^\prime(r)}{-V^\prime(r)}=\frac{-(\sinh r)^{N-1}v_n^\prime(r)}{-(\sinh r)^{N-1}V^\prime(r)} \ge \frac{C_V (1-\epsilon) }{C_V} \ge 1-\epsilon \ \ \ \forall r \in (r_\epsilon,n) \, , \ \forall n \ge n(r_\epsilon,\epsilon) \, ,
\end{equation*}
that is \eqref{eq: statProbNDeriv} for the sequence $\{ v_n \}$.

In general it is not guaranteed that $v_n(r) \le V(r)$ for all $r\in [1,\infty)$. Therefore, in order to conclude our proof, it is necessary to modify $\{ v_n \}$. To this aim, consider the following sequence:
$$ V_n=\lambda_n \, v_n \, , $$
where, for any fixed $n$, $\lambda_n$ is the \emph{largest} number belonging to $(0,1]$ for which $\lambda_n \, v_n(r) \le V(r)$ for all $r\in [1,\infty)$. We can assume that $\lambda_n < 1$ eventually, otherwise there is nothing to prove since along a subsequence $\{V_n\}$ has all the properties claimed in the statement of the Lemma. But if $\lambda_n$ is strictly smaller than $1$ then $V_n(r)$ necessarily touches $V(r)$ at some point $r=\xi_n \in [1,\infty)$ (see figures \ref{fig: csi1} and \ref{fig: csiInt}), because otherwise $\lambda_n$ would not be the largest number in $(0,1)$ for which $\lambda_n \, v_n(r) \le V(r)$ for all $r\in [1,\infty)$ (recall that each $v_n$ is compactly supported, $V$ is strictly positive and both $ v_n $ and $V$ are continuous). Now there are two possibilities: either the sequence $\{ \xi_n \}$ remains bounded or it is unbounded. In the first case, along subsequences $\{\xi_n\}$ converges to a certain value $\bar{\xi} \in [1,\infty)$. Since $\{v_n\}$ also converges \emph{locally uniformly} in $[1,\infty)$ to $V$ (see Remark \ref{oss: locUnif} below), then
$$   v_n(\xi_n)  \rightarrow V(\bar{\xi}) >0 \, ; $$
but by definition of $\xi_n$,
$$ \lambda_n\,v_n(\xi_n)   =   V(\xi_n)  \rightarrow V(\bar{\xi}) \, .  $$
Hence $ \lambda_n \rightarrow 1$. In the case $ \xi_n \rightarrow \infty$ (again, along subsequences), clearly each $\xi_n$ lies in the interior of $[1,n)$ eventually. Therefore, in addition to $V_n(\xi_n)=V(\xi_n)$, we also have $V_n^\prime(\xi_n)=V^\prime(\xi_n)$. So,
$$ - \lambda_n \, (\sinh \xi_n)^{N-1} \, v_n^\prime(\xi_n) = -(\sinh \xi_n)^{N-1} \, V_n^\prime(\xi_n) = -(\sinh \xi_n)^{N-1} \, V^\prime(\xi_n) \rightarrow C_V > 0 \, , $$
but also
$$ -(\sinh \xi_n)^{N-1} \, v_n^\prime(\xi_n) = {\int}_0^{\xi_n} \frac{v_n^{\frac{1}{m}}(s)}{T_n} \, (\sinh s)^{N-1} \mathrm{d}s \rightarrow \int_0^\infty \, \frac{V^{\frac{1}{m}}(s)}{T} \, (\sinh s)^{N-1} \mathrm{d}s  = C_V \, , $$
\begin{figure}
\includegraphics[scale=0.5]{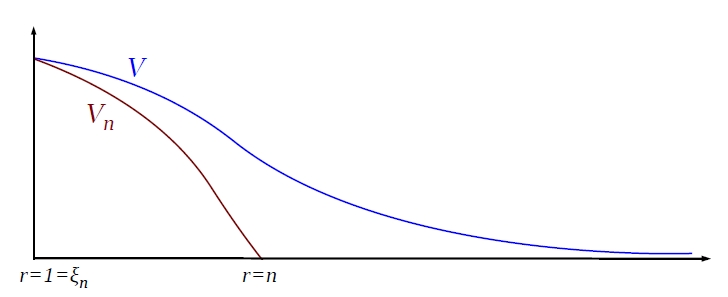}
\caption{$V_n$ touches $V$ at $\xi_n=1$.}
\label{fig: csi1}
\end{figure}
\begin{figure}
\includegraphics[scale=0.5]{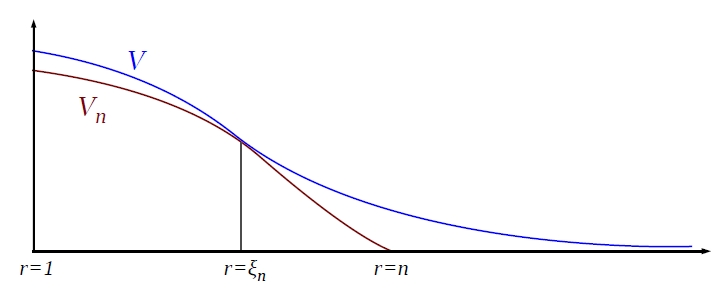}
\caption{$V_n$ touches $V$ at $\xi_n>1$.}
\label{fig: csiInt}
\end{figure}

\noindent where we have used once again the fact that $\{v_n\}$ converges to $V$ in $L^{{1}/{m}}_\textnormal{rad}$ and $ T_n \uparrow T$. This means that in this case as well $\{\lambda_n\}$ is forced to go to $1$ and so $\{V_n\}$ is indeed a sequence which has all the required properties (no subsequence is needed since the ongoing proof actually holds along any subsequence). Just observe that the parameter $T_n$ appearing in \eqref{eq: statProbN} here is actually $ T_n \, \lambda_n^{{1}/{m}-1}$.
\end{proof}
\end{lem}
\begin{oss}\label{oss: locUnif} \rm
As anticipated in the proof above, convergence of the sequence $\{v_n\}$ solving \eqref{eq: probMin} to $V$ also occurs \emph{locally uniformly} in $[1,\infty)$. This fact can be proved in several ways. One of them is the following: $\{v_n\}$ converges weakly in $\Hr$ to $V$ (so strongly in $L^{(m+1)/m}_\textnormal{rad}$), but using the equation solved by each $v_n$ and the fact that $ T_n \uparrow T$ one gets $ \| v_n^\prime\|_{2} \rightarrow \| V^\prime \|_{2}$. Therefore such convergence is actually strong. Since $\Hr$ is locally (in $(0,\infty)$) continuously embedded in $L^\infty$, the assertion holds.
\end{oss}
For reasons that will become clearer later (see the proof of Proposition \ref{thm: LimSup} below), instead of analysing minus the relative error $\psi$ it is convenient to study its natural approximation $\psi_n=1-{w^m}/{V_n}$, where $\{V_n\}$ is the sequence coming from Lemma \ref{lem: Approx}. It is straightforward to check that $\psi \ge \psi_n$ and the equation solved by $\psi_n$ is the following:
\begin{equation}\label{eq: MenErrRelN}
\frac{1}{m} \left(1-\psi_n \right)^{\frac{1}{m}-1} \psi_{n,\tau} = V_n^{1-\frac{1}{m}}\Delta \psi_n + 2 \, \frac{V_n^\prime \psi_n^\prime}{V_n^{\frac{1}{m}}} + \frac{1}{(1-m)} \left[\frac{1}{T_n} \left(1-\psi_n \right) - \frac{1}{T} \left(1-\psi_n \right)^{\frac{1}{m}} \right] ,  % \ \ \ \forall(r,t) \in (0,\infty) \times (1,\infty)
\end{equation}
valid in the region $\{(r,\tau) \in (0,n) \times (0,\infty)\}$.

The following fundamental Lemma provides a family of supersolutions to \eqref{eq: MenErrRelN}.
\begin{lem}\label{thm: MenErrRelNBarr}
Let $ \tau_0 >0 $. Consider the function
\begin{equation}\label{eq: formaPsi}
\Psi(r,\tau)=C-\frac{B}{r}-A\left(\tau-\tau_0\right) .
\end{equation}
If the positive parameters $A$, $B$, $C$, $\bar{r}$ satisfy the condition
\begin{equation}\label{eq: PsiFinalCondPre}
B K_\epsilon \, e^{(N-1)(1-\epsilon)\left(\frac{1}{m}-1 \right) \bar{r}} \ge \frac{A}{m} + \frac{2}{(1-m)T}  \, ,
\end{equation}
where $\epsilon>0$ is a small fixed number ($\epsilon=1/5$ will work) and $K_\epsilon$ is a positive constant, then $\Psi$ is a supersolution to \eqref{eq: MenErrRelN} in the region $\{(r,\tau): \, \Psi(r,\tau) \in [0,1) \} \cap \{(r,\tau) \in (\hat{r},n) \times (\tau_0,\infty)\}$ independently of $n$ provided
\begin{equation*}\label{eq: condRN}
\hat{r} \ge \bar{r} \vee r^\ast  \, , \ \ n \ge n^\ast \, ,
\end{equation*}
$r^\ast$ and $n^\ast$ being positive numbers depending only on $\epsilon$, $m$ and $N$, where we use the notation $a\vee b:=\max(a,b)$.
\begin{proof}
Let us compute the derivatives of $\Psi$:
\begin{equation*}\label{eq: derBarrN}
\begin{aligned}
\Psi_\tau &= -A  \, , \ \Psi^\prime= \frac{B}{r^2} \, , \\
\Delta \Psi &= -B\Delta\left( \frac{1}{r} \right) = -B \left( \frac{2}{r^3} - (N-1) \frac{ \coth(r)}{r^2} \right) .
\end{aligned}
\end{equation*}
In order to ensure that $\Psi$ is a supersolution to \eqref{eq: MenErrRelN} it is enough to ask
$$ \frac{1}{m} \left(1-\Psi \right)^{\frac{1}{m}-1} \Psi_{\tau} \ge V_n^{1-\frac{1}{m}}\Delta \Psi + 2 \, \frac{V_n^\prime \Psi^\prime}{V_n^{\frac{1}{m}}} + \frac{1}{(1-m)T_n}  \left(1-\Psi \right) \, , $$
that is
\begin{equation}\label{derBarrN1}
-\frac{A}{m} \left(1-\Psi \right)^{\frac{1}{m}-1}  \ge \underbrace{-\frac{B}{V_n^\frac{1}{m}} \left[ V_n \left( \frac{2}{r^3} - (N-1)  \frac{\coth(r)}{r^2} \right) -2 \frac{V_n^\prime}{r^2}  \right] + \frac{1}{(1-m)T_n}  \left(1-\Psi \right) }_{\mathcal{R}} \, ,
\end{equation}
where we have dropped the nonpositive term $-(1-\Psi)^{{1}/{m}}$. Now we need to suitably estimate from above the r.h.s.\ of \eqref{derBarrN1}. To this end, take $\epsilon>0$ sufficiently small. Moreover, take $n_0$ and $r_0$ so large that $T_n\ge T/2$ for all $n \ge n_0$ and $\coth(r) \le 1+ \epsilon$ for all $ r \ge r_0 $. Thanks to \eqref{eq: statProbNDeriv} and to the fact that $V_n \le V$ in $[1,\infty)$ we have:
\begin{equation*}\label{derBarrN2}
\begin{aligned}
\mathcal{R} & \le \frac{B}{V_n^\frac{1}{m}} \left[ (N-1)(1+\epsilon) \frac{V_n}{r^2} +2 \frac{V_n^\prime}{r^2}  \right] + \frac{2}{(1-m)T}  \left(1-\Psi \right) \\
 & \le \frac{B}{r^2 V_n^\frac{1}{m}} \left[ (N-1)(1+\epsilon) V +2 (1-\epsilon) V^\prime \right] + \frac{2}{(1-m)T}  \left(1-\Psi \right)
\end{aligned}
\end{equation*}
$$  \forall (r,\tau) \in (r_0\vee r_\epsilon \vee 1,n) \times (\tau_0,\infty) \, , \ \ \  \forall  n \ge n_0 \vee n_\epsilon \, . $$
Recall that the stationary solution $V$ satisfies \eqref{eq: fundEqIntLim}; in particular, by \eqref{eq:asympt}, there exists $r_1$ sufficiently large such that
\begin{equation*}\label{derBarrU}
V^\prime(r) \le -(N-1)(1-\epsilon) V(r) \ \ \ \forall r \ge r_1 \, .
\end{equation*}
Therefore,
\begin{equation}\label{derBarrN3}
\mathcal{R} \le -(N-1)\frac{B V}{r^2 V_n^\frac{1}{m}} \left[ -(1+\epsilon) + 2 (1-\epsilon)^2  \right] + \frac{2}{(1-m)T}  \left(1-\Psi \right)
\end{equation}
$$  \forall (r,\tau) \in (r_0\vee r_1 \vee r_\epsilon \vee 1,n) \times (\tau_0,\infty) \, , \ \ \  \forall  n \ge n_0 \vee n_\epsilon \, . $$
Since $\epsilon$ is small (it is enough that $ -(1+\epsilon) + 2 (1-\epsilon)^2 > 0$) and $V_n \le V$ in $[1,\infty)$, we can replace $V_n$ by $V$ in \eqref{derBarrN3}, which yields
\begin{gather}\label{derBarrN4}
\mathcal{R} \le -\frac{B}{r^2 V^{\frac{1}{m}-1}} \underbrace{(N-1) \left[ -(1+\epsilon) + 2 (1-\epsilon)^2  \right]}_{K_\epsilon} + \frac{2}{(1-m)T}  \left(1-\Psi \right) \\
\forall (r,\tau) \in (r_0\vee r_1 \vee r_\epsilon \vee 1 ,n) \times (\tau_0,\infty) \, , \ \ \  \forall  n \ge n_0 \vee n_\epsilon  \, . \nonumber
\end{gather}
Collecting \eqref{derBarrN1} and \eqref{derBarrN4}, we deduce that a necessary condition for $\Psi$ to be a supersolution to \eqref{eq: MenErrRelN} in the region $\{(r,\tau): \, \Psi(r,\tau) \in [0,1) \} \cap \{(r,\tau) \in (r_0\vee r_1 \vee r_\epsilon \vee 1 ,n) \times (\tau_0,\infty) \}$ is
\begin{equation}\label{derBarrN5}
\frac{B}{r^2 V^{\frac{1}{m}-1}} K_\epsilon  \ge \frac{A}{m} \left(1-\Psi \right)^{\frac{1}{m}-1} + \frac{2}{(1-m)T}  \left(1-\Psi \right) \, .
\end{equation}
Still from \eqref{eq: fundEqIntLim} and \eqref{eq:asympt} we get that
$$ \frac{K_\epsilon}{r^2 V^{\frac{1}{m}-1}} \ge K_\epsilon \, e^{(N-1)(1-\epsilon)\left(\frac{1}{m}-1 \right) r}  \ \ \ \forall r \ge r_2  $$
for another positive constant $K_\epsilon$, as long as $r_2$ is large enough. So, the final condition on $\bar{r}$ is
\begin{equation}\label{eq: PsiFinalCond}
B K_\epsilon \, e^{(N-1)(1-\epsilon)\left(\frac{1}{m}-1 \right) \bar{r}} \ge \frac{A}{m} + \frac{2}{(1-m)T}  \, .
\end{equation}
%provided $\bar{r} \ge r_0 \vee r_1 \vee r_2 \vee r_\epsilon \vee 1 $.
Summing up, if one fixes the positive parameters $\epsilon$ (small enough, say $\epsilon=1/5$), $A$, $B$, $C$, chooses ${r}^\ast = r_0 \vee r_1 \vee r_2 \vee r_\epsilon \vee 1$ and takes $ \bar{r}$ so large that \eqref{eq: PsiFinalCond} holds, then the function $\Psi$ as in \eqref{eq: formaPsi} is a supersolution to \eqref{eq: MenErrRelN} in the region $\{(r,\tau): \, \Psi(r,\tau) \in [0,1) \} \cap \{(r,\tau) \in (\hat{r},n) \times (\tau_0,\infty) \} $ for all $\hat{r} \ge \bar{r} \vee r^\ast$ and $n \ge n_0 \vee n_\epsilon = n^\ast $.
\end{proof}
%${r}^\ast = r_0(\epsilon) \vee r_1(\epsilon,m,N) \vee r_2(\epsilon,m,N) \vee r_\epsilon(m,N) \vee 1$
\end{lem}
In the next theorem we shall see, along the lines of the proof of \cite[Th. 2.1]{BGVFast}, how to use the (family of) barriers provided by Lemma \ref{thm: MenErrRelNBarr} in order to  prove that $\psi$ becomes smaller and smaller as $\tau \rightarrow \infty$. To this aim it is crucial to recall that Theorem \ref{CorHarnack}, which we proved in Sections \ref{Sec: Above}--\ref{sec: below}, implies the existence of a time $\tau_{w}>0$ and two positive constants $c_0,c_1$ such that
\begin{equation}\label{eq: HarnackPrinc}
c_0 \le \frac{w^m(r,\tau)}{V(r)} \le c_1 \ \ \ \forall(r,\tau) \in [0,\infty) \times [\tau_w, \infty) \, , \ \ c_0 < 1 < c_1 \, .
\end{equation}
\begin{pro}\label{thm: LimSup}
Let $\psi=1-{w^m}/{V}$. There holds
\begin{equation*}
\limsup_{\tau \rightarrow \infty}  \sup_{r \in [0,\infty)} \psi(r,\tau) \le 0   \, .
\end{equation*}
\begin{proof}
First of all consider the barrier $\Psi$ as in \eqref{eq: formaPsi}, with the choices $A=B=1$. Then set $C=1-{c_0}/{2}$ and assume that $\hat{r}$ is greater than $(2/c_0) \vee \bar{r} \vee r^\ast$, where $\bar{r}$ is taken large enough so as to satisfy
\begin{equation*}
K_\epsilon \, e^{(N-1)(1-\epsilon)\left(\frac{1}{m}-1 \right) \bar{r}} \ge \frac{1}{m} + \frac{2}{(1-m)T}  \, ,
\end{equation*}
that is \eqref{eq: PsiFinalCondPre} when $ A=B=1 $. By construction, such a $\Psi$ is always lower than $1$. Therefore, from Lemma \ref{thm: MenErrRelNBarr}, it is a supersolution to \eqref{eq: MenErrRelN} in the region $(\hat{r},n) \times (\tau_0,\infty)$ (let $\tau_0 \ge \tau_w$) for all $n \ge n^\ast$ as long as it is greater than or equal to zero. Since we want to compare $\Psi$ with the solutions $\psi_n$, let us check conditions on a parabolic boundary of the form $ [\hat{r},n] \times \left\{ \tau_0 \right\} \cup \left\{ \hat{r} \right\} \times (\tau_0,\tau_1) \cup \left\{ n \right\} \times (\tau_0,\tau_1) $, where $\tau_1>\tau_0$ is a positive time to be chosen later. On the bottom we have that, thanks to the choice of $C$ and \eqref{eq: HarnackPrinc},
\begin{equation*}\label{eq: condBott}
\Psi(r,\tau_0) = C -\frac{1}{r} \ge 1-c_0 \ge \psi(r,\tau_0) \ge \psi_n(r,\tau_0) \ \ \ \forall r \in [\hat{r},n] \, ,
\end{equation*}
the last inequality following from the fact that $V \ge V_n$ for $r\ge 1$. On the inner lateral boundary there holds
\begin{equation*}\label{eq: innBound}
\Psi(\hat{r},\tau) = C -\frac{1}{\hat{r}} - (\tau-\tau_0) = 1-\frac{c_0}{2} -\frac{1}{\hat{r}} - (\tau-\tau_0) \ \ \ \forall \tau \in (\tau_0,\tau_1) \, .
\end{equation*}
Now let us fix a small $\varepsilon > 0$ (let $ \varepsilon \le c_0/2$). Assume that $\hat{r}$ is also larger than $2/\varepsilon$. Exploiting the \emph{local} uniform convergence of $ \psi $ to zero (Proposition \ref{teo: convLocZero}), we know that if $\tau_0 \ge \tau_0(\hat{r},\varepsilon)$ then
$$  \psi(\hat{r},\tau) \le \frac{\varepsilon}{2} \ \ \ \forall \tau \ge \tau_0 \, .  $$
Hence, once $\tau$ satisfies
$$ \tau \le \tau_1 = \tau_0 + 1 - \frac{c_0}{2}-\frac{1}{\hat{r}}-\frac{\varepsilon}{2} \ge \tau_0 + 1-c_0 > \tau_0 \, , $$
it is guaranteed that
\begin{equation*}\label{eq: latBound}
\Psi(\hat{r},\tau) \ge \Psi(\hat{r},\tau_1) = \frac{\varepsilon}{2} \ge \psi(\hat{r},\tau) \ge \psi_n(\hat{r},\tau)  \ \ \ \forall \tau \in (\tau_0,\tau_1) \, .
\end{equation*}
On the outer lateral boundary $\{ n \} \times (\tau_0,\tau_1)$ it is plain that $\Psi(n,\tau) > \psi_n(n,\tau)=-\infty $. In fact this is the reason why we needed to suitably approximate $V$ from below with the sequence $V_n$, otherwise we would have not been able to compare $\Psi$ with $\psi$ in an outer lateral boundary. So, from the comparison principle%\footnote{By construction, $\Psi(r,t) \ge \frac{\epsilon}{2} \ge 0$ for all $(r,t) \in [\bar{r},n]\times[t_0,t_1] $.}
, we get that $\Psi \ge \psi_n $ in the region $[\hat{r},n]\times[\tau_0,\tau_1] $. In particular,
\begin{equation}\label{eq: AfterCompN}
\psi_n(r,\tau_1) \le \Psi(r,\tau_1) \le \varepsilon    \ \ \ \forall r \in [\hat{r},n] \, .
\end{equation}
Passing to the limit in \eqref{eq: AfterCompN} as $n \rightarrow \infty$ yields
\begin{equation}\label{eq: AfterComp}
\psi(r,\tau_1) \le \Psi(r,\tau_1) \le \varepsilon    \ \ \ \forall r \in [\hat{r},\infty) \, .
\end{equation}
Since $\tau_1 = \tau_0 + 1 - {c_0}/{2}-{1}/{\hat{r}}-{\varepsilon}/{2}$ and \eqref{eq: AfterComp} holds for all $\tau_0 \ge \tau_0(\hat{r},\varepsilon)$, we have:
\begin{equation}\label{eq: AfterCompBis}
\psi(r,\tau) \le \varepsilon  \ \ \ \forall (r,\tau) \in [\hat{r},\infty) \times \left[\bar{\tau}, \infty \right) ,
\end{equation}
where $\bar{\tau} = \tau_0(\hat{r},\varepsilon) + 1 - {c_0}/{2}-{1}/{\hat{r}}-{\varepsilon}/{2}$. Thanks to \eqref{eq: AfterCompBis} and to the local uniform convergence of $\psi$ to zero, we conclude that
\begin{equation}\label{eq: LimSupEps}
\limsup_{\tau \rightarrow \infty}  \sup_{r \in [0,\infty)} \psi(r,\tau) \le \varepsilon \, .
\end{equation}
Because \eqref{eq: LimSupEps} is true for all $\varepsilon>0$, it holds for $\varepsilon=0$ too and the proof is complete.
\end{proof}
\end{pro}
The final step is to prove an analogous result for the relative error $\phi={w^m}/{V}-1$. Since the ideas are similar to the ones developed in the proofs of Lemma \ref{thm: MenErrRelNBarr} and Proposition \ref{thm: LimSup}, we shall proceed in a concise way.
\begin{lem}\label{thm: PlusErrRelBarr}
Let $ \tau_0 > 0 $. Consider the function
\begin{equation}\label{eq: formaPhi}
\Phi(r,\tau)=C-\frac{B}{r}-A\left(\tau-\tau_0\right) .
\end{equation}
If the positive parameters $A$, $B$, $C$, $\bar{r}$ satisfy the condition
\begin{equation}\label{eq: PhiFinalCondPre}
B K_\epsilon \, e^{(N-1)(1-\epsilon)\left(\frac{1}{m}-1 \right) \bar{r}} \ge \left(1 + C \right)^{\frac{1}{m}-1} \left( \frac{A}{m} + \frac{2(1+C)}{(1-m)T} \right) ,
\end{equation}
where $\epsilon>0$ is a small fixed number ($\epsilon=1/5$ will work) and $K_\epsilon$ is a positive constant, then $\Phi$ is a supersolution to \eqref{eq: PiuErrRel} in the region $\{(r,\tau): \, \Phi(r,\tau) > -1 \} \cap \{(r,\tau) \in (\hat{r},\infty) \times (\tau_0,\infty)\}$, provided $\hat{r} \ge  \bar{r} \vee {r}_\ast$, where ${r}_\ast={r}_\ast(\epsilon,m,N) > 0$.
\begin{proof}
Condition \eqref{eq: PhiFinalCondPre} can be obtained reasoning as in the proof of Lemma \ref{thm: MenErrRelNBarr}. Indeed it is even easier to get to an inequality like \eqref{derBarrN5} with $(1-\Psi)^{1/m-1}$ replaced by $(1+\Phi)^{1/m-1}$ and $(1-\Psi)$ replaced by $(1+\Phi)^{1/m}$, because here we need not deal with the approximating sequence $V_n$. Then, in order to arrive at \eqref{eq: PhiFinalCondPre}, one just uses the fact that by construction $\Phi \le C$.
%The requirement $ \Phi(r,t) > -1 $ is equivalent to $ \Psi(r,t) < 1 $ in Lemma \ref{thm: MenErrRelNBarr} : outside this region the equation for which $\Phi$ is a supersolution has the wrong sign.
\end{proof}
\end{lem}
Thanks to the (family of) barriers \eqref{eq: formaPhi}, we can now prove the analogue of Proposition \ref{thm: LimSup} for $\phi$.
\begin{pro}\label{thm: LimSupPhi}
Let $\phi={w^m}/{V}-1$. There holds
\begin{equation*}
\limsup_{\tau \rightarrow \infty}  \sup_{r \in [0,\infty)} \phi(r,\tau) \le 0   \, .
\end{equation*}
\begin{proof}
Again, we proceed along the lines of the proof of Proposition \ref{thm: LimSup}. Fix $A=B=1$, $C=c_1-1/2$, $\varepsilon < 1/2$ and $\hat{r}$ large enough so that it satisfies $\hat{r} \ge  \bar{r} \vee {r}_\ast \vee 2/\varepsilon $, where $\bar{r}$ complies with \eqref{eq: PhiFinalCondPre}. These choices ensure that $\Phi$ is a supersolution to \eqref{eq: PiuErrRel} in the region $(\hat{r},\infty) \times (\tau_0,\tau_1)$, with
$$ \tau_1 = \tau_0 + c_1 - \frac{1}{2}-\frac{1}{\hat{r}}-\frac{\varepsilon}{2} \, . $$
Then, one takes $\tau_0 = \tau_0(\hat{r},\varepsilon)$ so that $\phi(\hat{r},\tau) \le \varepsilon/2 $ for all $\tau \ge \tau_0(\hat{r},\varepsilon)$; in this way $\Phi$ and $\phi$ are correctly ordered on $[\hat{r},\infty ) \times \{ \tau_0 \} \cup \{ \hat{r} \} \times (\tau_0,\tau_1)$. However, we have no control over their relation in an outer lateral boundary of the form $\{ n \} \times (\tau_0,\tau_1)$, for $n > \hat{r}$. In order to overcome this difficulty it is enough to replace $w$ by $w_n \le w$, the latter being the sequence of rescaled solutions to \eqref{eq: Fastn}, which vanish on $\{ n \} \times (\tau_0,\tau_1)$. By construction $\phi_n=w_n^m/V-1 \le \phi$ and the equation solved by $\phi_n$ is basically the same as \eqref{eq: PiuErrRel} (just replace ${(1+\phi)^{1/m}}/{T}$ with ${(1+\phi_n)^{1/m}}/{T_n}$), so that for $n$ large (that is, $T_n$ close to $T$) $\Phi$ is a supersolution also to such equation. Therefore from now on one proceeds exactly as in the proof of Proposition \ref{thm: LimSup}.
\end{proof}
\end{pro}
Going back to the original variables $t$ and $u(r,t)$, we have that Propositions \ref{thm: LimSup} and \ref{thm: LimSupPhi} give convergence of $u^m(r,t)$ to ${(1-t/T)^{m/(1-m)}}V(r)$ in relative error. From this fact the first main result stated in Theorem \ref{ThConvRel}, namely formula \eqref{eq: ConvRelFin}, clearly follows.

\section{Results for derivatives}\label{Sec: deriv}
This section is devoted to proving the claimed results of Theorems \ref{ThConvRel}, \ref{CorHarnack} and \ref{CorSub} which deal with derivatives. In order to do that, we shall exploit a useful change of spatial variable and rescaling techniques in the spirit of \cite{DKV91}.

Let
\begin{equation}\label{eq: StoR}
s=s(r)=\int_r^{+\infty} \frac{1}{(\sinh \omega)^{N-1}} \, \mathrm{d}\omega \, ,
\end{equation}
so that \[\frac{\mathrm{d}}{\mathrm{d}s}=-(\sinh r)^{N-1} \frac{\mathrm{d}}{\mathrm{d}r}.\]
Notice that $ s \downarrow 0 $ as $r \uparrow \infty $, while $s \uparrow \infty $ as $r \downarrow 0 $.
%(recall that $N \ge 2$ and $\sinh(\omega) \asymp \omega$ as $\omega \rightarrow 0$).
Hence, the hyperbolic Laplacian of a regular function $f=f(r)=f(r(s))$ is
\begin{equation*}\label{eq: hypLapS2}
\Delta f = (\sinh r)^{-2(N-1)} \, \frac{\mathrm{d}^2 f}{\mathrm{d}s^2 }  \, .
\end{equation*}
Let $w$ be the solution to \eqref{eq: FastRisc}. Take $\tau_0\ge 1$ and $r_0\ge R$ with $R$ large enough (to be chosen later), and let $s_0=s(r_0)$. Consider the rescaled function
%(about $(r_0,\tau_0)$)
\begin{equation}\label{eq: hypLapS}
\mathcal{W}(x,y)= (\sinh r_0)^{\frac{N-1}{m}}\, w \! \left(r(s_0 + x \, (\sinh r_0)^{-(N-1)}),\tau_0 + y \, (\sinh r_0)^{-(N-1)\left(\frac1m-1\right)}\right) ,
\end{equation}
where $(x,y)$ varies in the square $(-\alpha,\alpha)^2$, $\alpha>0$ being a fixed number to be set later. After straightforward computations one sees that $\mathcal{W}$ satisfies the following equation:
\begin{equation}\label{eq: riscRisc}
\mathcal{W}_y = \underbrace{ \frac{(\sinh r_0)^{2(N-1)}}{(\sinh r)^{2(N-1)}} }_{a(r,r_0)} \left( \mathcal{W}^m \right)_{xx} + \underbrace{ \frac{(\sinh r_0)^{-(N-1)\left(\frac1m-1\right)}}{(1-m)T} }_{b(r_0)} \, \mathcal{W}  \, .
\end{equation}
From now on we shall fix $\alpha=\alpha(N)$, independent of $r_0$ and $\tau_0$, so small that $s$ is forced not to go out of the interval $(s_0/2,3s_0/2)$ as $x$ varies in the interval $(-\alpha,\alpha)$. Such a choice of $\alpha$ is feasible. Indeed, notice first that from \eqref{eq: StoR} one has
\begin{equation}\label{eq: limSR1}
\lim_{r \rightarrow \infty} \frac{s(r)}{(\sinh r)^{-(N-1)}} = C(N)>0 \, .
\end{equation}
Moreover, $s$ varies in the interval $ (s_0 - \alpha (\sinh r_0)^{-(N-1)},s_0 + \alpha (\sinh r_0)^{-(N-1)}) $. Hence, if $r_0 \ge R$ (namely, $s_0 \le s(R)$) with $R=R(N)$ large enough, a proper choice of $\alpha=\alpha(N)$ (sufficiently small) ensures that $s$ does not leave the interval $(s_0/2,3s_0/2)$ as $x \in (-\alpha,\alpha)$. Still from \eqref{eq: limSR1} one deduces that
\begin{equation}\label{eq: limSR2}
\lim_{s \rightarrow 0} r(s/2) - r(s) = \underline{C} \, , \ \ \ \lim_{s \rightarrow 0} r(3s/2) - r(s) = \overline{C}
\end{equation}
for suitable constants $\underline{C}$ and $\overline{C}$. The choice of $\alpha$ above and \eqref{eq: limSR2} imply that, for $r_0=r(s_0)$ large enough (again, greater or equal than a given value $R$ that depends only on $N$), there holds
\begin{equation}\label{eq: boundsSR}
 \frac{1}{C_1} \sinh r_0  \le \sinh r \le C_1 \sinh r_0  \ \ \  \forall r \in \left( r(s_0/2) \, , r(3s_0/2) \right) ,
\end{equation}
$ C_1 $ being another positive constant. Gathering all this information, let us go back to equation \eqref{eq: riscRisc}. First, from the global Harnack principle \eqref{eq: ghp} and from \eqref{eq: boundsSR}, one deduces that
\begin{equation}\label{eq: boundWrisc}
| \mathcal{W}(x,y) | \le M \ \ \ \forall (x,y) \in (-\alpha,\alpha)^2
\end{equation}
for a positive constant $M=M(u_0,m,N)$ that does not depend on $r_0,\tau_0$. Moreover, notice that
\begin{equation}\label{eq: derXR}
\frac{\mathrm{d}}{\mathrm{d}x} = - \frac{(\sinh r)^{N-1}}{(\sinh r_0)^{N-1}}  \frac{\mathrm{d}}{\mathrm{d}r} \, .
\end{equation}
%That is, derivating w.r.t. $x$ is basically the same as derivating w.r.t. $r$ as long as $r$ varies in a small neighbourhood of $r_0$.
Hence, \eqref{eq: riscRisc}, \eqref{eq: derXR} and the inequality
$$ \frac{1}{C_1} \cosh r_0  \le \cosh r \le C_1 \cosh r_0 \ \ \ \forall r \in \left( r(s_0/2), r(3s_0/2) \right) , $$
which can be proved as \eqref{eq: boundsSR} above, give
\begin{equation}\label{eq: boundWriscCoeff}
\frac{1}{P_k} \le  \left| \frac{\mathrm{d}^k a(r(x),r_0) }{\mathrm{d}x^k} \right| \le P_k \, , \ \ \left| b(r_0) \right| \le Q  \ \ \ \ \forall (x,y) \in (-\alpha,\alpha)^2 \, , \ \forall k \in \mathbb{N} \cup \{ 0 \} \, ,
\end{equation}
again for positive constants $P_k=P_k(N)$ and $Q=Q(T(u_0),m,N)$ that do not depend on $r_0,\tau_0$. The bounds \eqref{eq: boundWrisc} and \eqref{eq: boundWriscCoeff} permit to conclude, as in the proofs of Theorem 1.1 and Lemmas 3.1-3.2 of \cite{DKV91} (one uses the regularity results of \cite{Fri} and \cite{CDib} for bounded solutions to nonlinear parabolic equations like \eqref{eq: riscRisc}), that the following estimates hold:
\begin{equation}\label{eq: stimeDerivateX}
\left| \frac{{\partial}^k \mathcal{W}}{\partial x^k}(0,0) \right|\le M_{1,k}  \, , \  \left| \frac{\partial^k \mathcal{W}}{\partial y^k}(0,0) \right|\le M_{2,k} \ \ \ \ \forall k \in \mathbb{N} \, ,
\end{equation}
$M_{1,k}$ and $M_{2,k}$ being suitable positive constants depending only on $M$, $\alpha$, $T$, $m$ and $N$. Going back to the original variables, we end up with the existence of positive constants $A_{1,k}=A_{1,k}(u_0,m,N,\varepsilon)$, $A_{2,k}=A_{2,k}(u_0,m,N,\varepsilon)$ such that
\begin{gather}
\left|\frac{\partial^k w}{\partial r^k}(r_0,\tau_0)\right| \le A_{1,k} \, e^{-\frac{N-1}{m}r_0} \, , \label{eq: stimeDerivateR1} \\
\left|\frac{\partial^k w}{\partial \tau^k}(r_0,\tau_0)\right| \le A_{2,k} \, e^{k(N-1)\left( \frac1m -1\right)r_0} \, e^{-\frac{N-1}{m}r_0}  \label{eq: stimeDerivateR2} \\
 \forall r_0 \ge 0 \, , \ \forall \tau_0\ge\varepsilon \, , \ \forall k \in \mathbb{N} \, . \nonumber
\end{gather}
In order to prove \eqref{eq: stimeDerivateR1}-\eqref{eq: stimeDerivateR2} in the region $r_0\ge R, \tau_0\ge1$ one just uses \eqref{eq: hypLapS}, \eqref{eq: derXR}, \eqref{eq: stimeDerivateX} and the equality $ \mathrm{d}\tau = (\sinh r_0)^{-(N-1)\left(\frac1m - 1 \right)} \mathrm{d}y$ recursively. To extend such bounds to $r_0 \in (0,R)$ it is enough to apply the aforementioned regularity results to $w$ directly, since in this region $w$ is bounded away from zero. Finally, they also hold for $\tau_0\ge \varepsilon$ up to admitting dependence of the constants $A_{1,k},A_{2,k}$ on $\varepsilon$ as well through $\alpha$ and $M$ (as a consequence of the global Harnack principle \eqref{eq: ghp}).

%Since $(r_0,\tau_0)$ is a generic point (one only needs to assume $r_0$ large enough and $\tau_0>0$), and estimates like \eqref{eq: stimeDerivateR} do hold \emph{locally} in space thanks to the uniform positivity of $w$ and so thanks to standard quasilinear theory ({\color{red} citare qualcosa}), we conclude that
%\begin{equation}\label{eq: stimeDerivateRFinale}
%\left| \frac{\mathrm{d}^k w}{\mathrm{d}r^k}(r,\tau) \right| \le D(k,\varepsilon) \, \sinh^{-\frac{N-1}{m}}(r) \ \ \ \forall (r,\tau) \in (0,\infty) \times (\varepsilon,\infty) \, , \ \forall k \in \mathbb{N} \, , \ \forall \varepsilon > 0 \, ,
%\end{equation}
%for positive constants $D=D(k,\varepsilon)$.
% or \eqref{eq: fundEqIntLim}), and exploiting recursively equation \eqref{eq: statProb} itself, we deduce that derivatives of $V(r)$ up to \emph{any} order behave like $ e^{-(N-1)r}$ at infinity (recall formula \eqref{eq:asympt}).
Using the fact that the solution $V(r)$ to \eqref{eq: statProb} and any of its derivatives $V^{(k)}(r)$ behave like $e^{-(N-1)r}$ at infinity (recall \eqref{eq:asympt} and see formula \eqref{eq: potenzeV} below), it is only a matter of tedious computations to check that estimates \eqref{eq: stimeDerivateR1} imply that
\begin{equation}\label{eq: stimeDerivateRelErr}
\left\| \frac{\partial^k \varphi}{\partial r^k}(\tau) \right\|_\infty = \left| \varphi(\tau) \right|_{C^k(\mathbb{R})} \le A_{1,k} \ \ \ \forall \tau \ge \varepsilon \, , \ \forall k \in \mathbb{N}
\end{equation}
(up to possibly different constants that we keep denoting as $A_{1,k}$), where $\varphi$ is the relative error $\varphi(r,\tau)={w(r,\tau)}/{V^{1/m}(r)}-1$ (notice that it is different from the relative error $\phi$ we dealt with in Section \ref{Sec: rel}). By exploiting the interpolation inequalities (we refer to \cite[p. 130]{Ni} or \cite[App. 3]{BGV2} for a short review on general interpolation inequalities)
$$ \left| \varphi(\tau) \right|_{C^k(\mathbb{R})} \le C_{j,k} \left| \varphi(\tau) \right|_{C^j(\mathbb{R})}^{k/j} \left\| \varphi(\tau) \right\|_\infty^{1-k/j} \ \ \ \forall k,j \in \mathbb{N}: \ k<j \, , $$
the bounds \eqref{eq: stimeDerivateRelErr} and \eqref{eq: ConvRelFin} (that is, the fact that $\| \varphi(\tau) \|_\infty \rightarrow 0$ as $\tau \rightarrow \infty$), we conclude that also
\begin{equation}\label{eq: errRelCk}
\lim_{\tau \rightarrow \infty} \left| \varphi(\tau) \right|_{C^k(\mathbb{R})} = 0 \ \ \ \forall k \in \mathbb{N} \, .
\end{equation}
Recalling \eqref{eq: rescSol}, this is equivalent to formula \eqref{deriv1} of Theorem \ref{ThConvRel}. As for estimates \eqref{eq: derivSP}-\eqref{eq: derivTE} of Theorem \ref{CorHarnack}, just notice that they can be readily deduced from \eqref{eq: stimeDerivateR1}-\eqref{eq: stimeDerivateR2} and again \eqref{eq: rescSol}.

Finally, let us investigate what \eqref{eq: errRelCk} means in terms of spatial derivatives for $w$ (or $u$). That is, we are going to prove formula \eqref{eq: convder} of Theorem \ref{ThConvRel}. First of all we show, by induction, that
\begin{equation}\label{eq: potenzeV}
\lim_{r \to \infty} \frac{\frac{\mathrm{d}^k V^\alpha}{\mathrm{d}r^k}(r) }{V^\alpha(r)} = \left(-\alpha (N-1)  \right)^k  \ \ \ \forall \alpha \neq 0 \, , \ \forall k \in \mathbb{N} \, .
\end{equation}
Indeed, the existence of the limit in \eqref{eq: potenzeV} is a consequence of the chain rule and \eqref{eq:asympt}. Then, suppose that \eqref{eq: potenzeV} holds for a given $k$. Consider the limit
\begin{equation}\label{eq: potenzeV0}
\lim_{r \to \infty} \frac{\frac{\mathrm{d}^{k+1} V^\alpha}{\mathrm{d}r^{k+1}}(r)} {\frac{\mathrm{d}V^\alpha}{\mathrm{d}r}(r) } = \frac1\alpha \lim_{r \to \infty} \frac{\frac{\mathrm{d}^{k+1} V^\alpha}{\mathrm{d}r^{k+1}}(r)}{V^\alpha(r)} \, \frac{V(r)}{V^\prime(r)} = \frac{1}{-\alpha(N-1)} \lim_{r \to \infty} \frac{\frac{\mathrm{d}^{k+1} V^\alpha}{\mathrm{d}r^{k+1}}(r)}{V^\alpha(r)}  \, .
\end{equation}
Since the limit on the r.h.s.\ of \eqref{eq: potenzeV0} exists, by de l'H\^{o}pital's Theorem it coincides with
$$ \lim_{r \to \infty} \frac{\frac{\mathrm{d}^k V^\alpha}{\mathrm{d}r^k}(r)}{V^\alpha(r)} \, ; $$
this proves the inductive step. Because \eqref{eq: potenzeV} is trivially valid for $k=0$, we conclude that it holds for all $k \in \mathbb{N}$.

% both $\| \varphi(\tau) \|_\infty$ and $| \varphi(\tau) |_{C^k(\mathbb{R})} $ go to zero as $\tau \rightarrow \infty$. Moreover,
Now recall that, for the $k$-th spatial derivative of $\varphi$, one has the binomial formula
\begin{equation}\label{eq: binomial}
\frac{\partial ^k \varphi}{\partial r^k} = \sum_{j=0}^{k} \binom{k}{j} w^{(j)} \left( V^{-\frac{1}{m}} \right)^{(k-j)} \, ,
\end{equation}
where, to simplify notation, we use apexes to denote derivatives w.r.t.\ $r$. We are going to show that, for all $k \in \mathbb{N}$, the ratio $w^{(k)}(r,\tau) / V^{1/m}(r,\tau)$ converges in $L^\infty$ to a smooth function $G_k(r)$ such that
\begin{equation}\label{eq: propGK}
\lim_{r \to \infty} G_k(r) = (-1)^k \left(\frac{N-1}{m} \right)^k \, .
\end{equation}
This is clearly equivalent to proving \eqref{eq: convder}. For $k=0$ \eqref{eq: propGK} is exactly \eqref{eq: ConvRelFin}, with the choice $G_0=1$. For further derivatives we shall prove that \eqref{eq: propGK} holds by induction. Indeed, suppose such result is true for all $j \le k-1$. From \eqref{eq: errRelCk} and \eqref{eq: binomial} we get
\begin{equation}\label{eq: binomial0}
\lim_{\tau \rightarrow \infty} \left\| \frac{w^{(k)}(\tau)}{V^{\frac1m}} - \sum_{j=0}^{k-1} -\binom{k}{j} w^{(j)}(\tau) \left( V^{-\frac{1}{m}} \right)^{(k-j)}   \right\|_\infty = 0 \, .
\end{equation}
Thanks to the inductive step and \eqref{eq: potenzeV}, \eqref{eq: binomial0} implies
\begin{equation*}\label{eq: binomial1}
\lim_{\tau \rightarrow \infty} \left\| \frac{w^{(k)}(\tau)}{V^{\frac1m}} - \sum_{j=0}^{k-1} -\binom{k}{j} G_j \frac{\left( V^{-\frac{1}{m}} \right)^{(k-j)}}{V^{-\frac{1}{m}} }   \right\|_\infty = 0 \, .
\end{equation*}
So, we are left with proving that the function
$$ G_k(r)= - \sum_{j=0}^{k-1} \binom{k}{j} G_j(r) \frac{\left( V^{-\frac{1}{m}} \right)^{(k-j)}(r)}{V^{-\frac{1}{m}}(r) } $$
complies with \eqref{eq: propGK}. From \eqref{eq: potenzeV} and the inductive hypothesis this will follow provided there holds
$$ -\left(\frac{N-1}{m} \right)^k \ \sum_{j=0}^{k-1} \binom{k}{j} (-1)^j  = (-1)^k \left(\frac{N-1}{m} \right)^k \, ,$$
an equality which is in fact true in view of the identity $\sum_{j=0}^{k-1} \binom{k}{j} (-1)^j  = (-1)^{k-1}$, valid by Newton's binomial formula.

\begin{oss}\label{oss: mSottocritico} \rm
In order to obtain the bounds \eqref{eq: derivSPCri}-\eqref{eq: derivTECri} of Theorem \ref{CorSub} (subcritical $m$), notice the quantities $c_1(t,\cdot)$, $c_2(t,\cdot)$ appearing in \eqref{eq: sub} can be taken to be continuous functions of $t \in (0,T)$, as a consequence of the method of proof of such inequality (see Sections \ref{Sec: Above} and \ref{sec: below}). This fact is sufficient to prove, as above, that there exists a constant $M$ (now depending on $\tau_0$ too) such that \eqref{eq: boundWrisc} holds true. This is enough in order to proceed along the previous lines and get \eqref{eq: derivSPCri}-\eqref{eq: derivTECri}.
\end{oss}

% binomial formula

%Note that \eqref{eq: ConvRelFin} is equivalent to
%\begin{equation}\label{eq: ConvRelFinDivm}
%\lim_{t \rightarrow \infty} \left\| \frac{u(r,t)}{\left( 1-\frac{t}{T} \right)^{\frac{1}{1-m}} V^{\frac{1}{m}}(r)} - 1 \right\|_\infty = 0 \, .
%\end{equation}

\bibliographystyle{plainnat}

%\addcontentsline{toc}{chapter}{Bibliografia}

\end{document}